\documentclass[british]{article}
\usepackage[T1]{fontenc}
\usepackage[utf8]{inputenc}
\usepackage{babel}
\usepackage{amsmath}
\usepackage{amssymb}
\usepackage{graphicx}
\usepackage{geometry}
\usepackage{microtype}
\usepackage[unicode=true,
 pdftitle={Geometric Functional Inequalities for Markov Processes: An Invitation via the Langevin Diffusion},
 pdfauthor={Sam Power},
 bookmarks=true,bookmarksnumbered=false,bookmarksopen=false,
 breaklinks=false,pdfborder={0 0 0},pdfborderstyle={},backref=false,colorlinks=false]
 {hyperref}

\date{}

\begin{document}
\title{Geometric Functional Inequalities for Markov Processes:\linebreak{}
An Invitation via the Langevin Diffusion}
\author{Sam Power\\
School of Mathematics, University of Bristol\\
\href{mailto:sam.power@bristol.ac.uk}{sam.power@bristol.ac.uk}}
\maketitle
\begin{abstract}
Functional inequalities provide a language for reasoning about the
long-time properties of ergodic Markov processes in terms of their
macroscopic geometry. This article offers an informal introduction
to functional inequalities for readers who are more familiar with
probabilistic tools based around couplings, Lyapunov techniques, and
similar. The main ideas are developed through the overdamped Langevin
diffusion, where various functionals related to this convergence admit
particularly transparent expressions, and for which the key associated
inequalities -- namely, the Poincaré and logarithmic Sobolev inequalities
-- are especially well-studied. We explain the qualitative content
of these inequalities, what they imply about convergence to equilibrium,
how they transform and compose under natural operations, and how they
interface with classical probabilistic approaches. Related inequalities
are discussed more briefly, and an annotated guide is provided for
readers wishing to enter the wider literature.
\end{abstract}
\begin{figure}[h]
\caption{A coarse impression of a Langevin diffusion in an asymmetric triple-well
potential.}

\centering{}\includegraphics[scale=0.25]{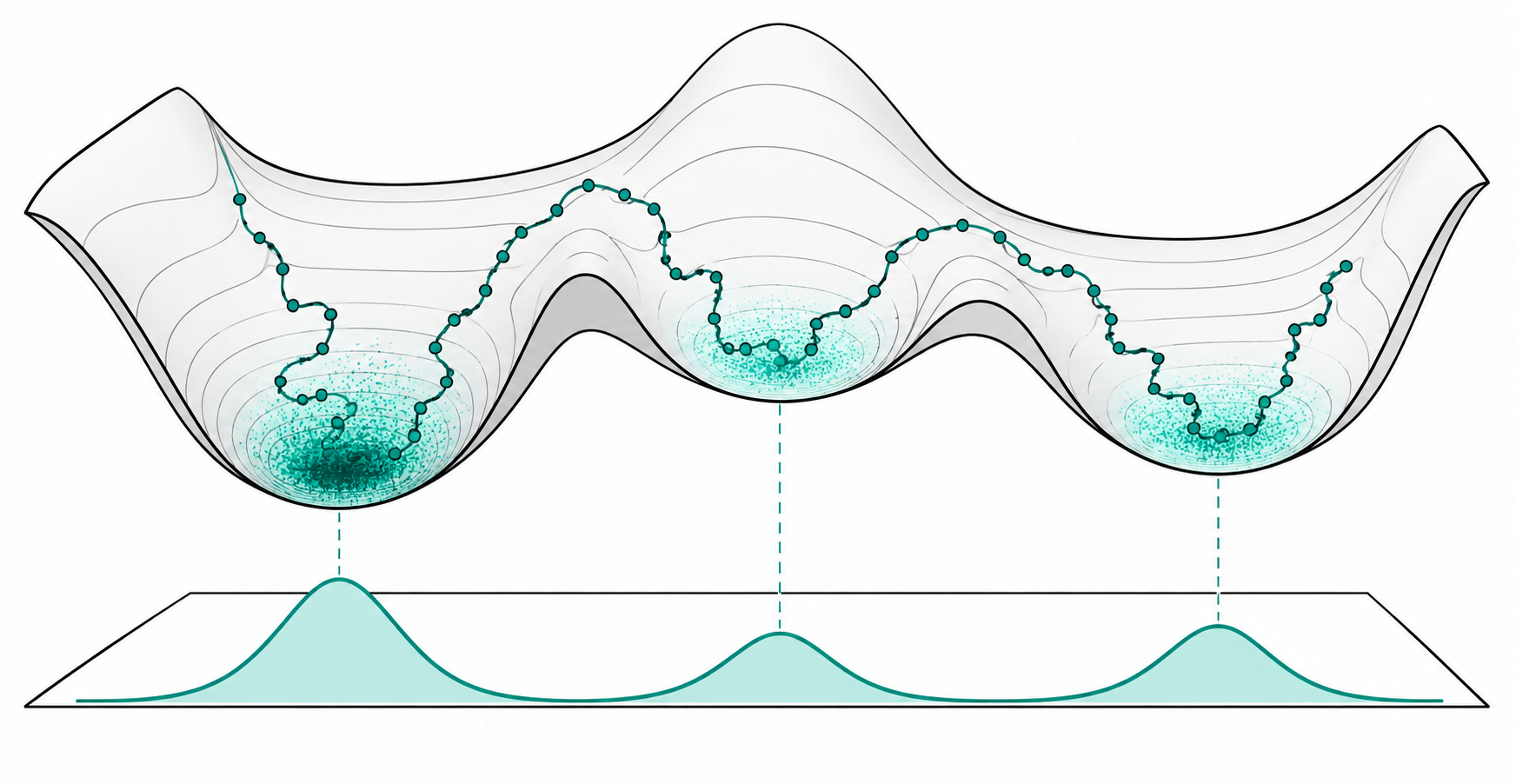}
\end{figure}

\newpage{}

\section{Prelude: Two Views of Long-Time Behaviour}
\begin{quote}
\begin{center}
So, how do \emph{you} like to study the long-time behaviour of a Markov
process?
\par\end{center}

\end{quote}
For many probabilists, a natural answer to this question is to study
the process in a pathwise manner: perhaps one constructs a coupling,
identifies regeneration times, proves a drift condition, or otherwise
determines what a typical trajectory of the process is doing. These
methods are often local and concrete, establishing them as useful
tools for structured problems in which this pathwise behaviour is
relatively transparent.

A functional-analytic approach begins at a different resolution. Instead
of following an individual trajectory, one steps back to study the
Markov semigroup directly, asking how it acts on functions or probability
measures. Convergence is then quantified through norms, divergences,
and entropies. Infinitesimal behaviour is encoded by a linear operator
on a suitable function space, and induces a bilinear energy form which
quantifies various pertinent notions of stationarity. The approach
is then relatively more global and coarse-grained by design, suppressing
many details of the path in order to expose large-scale structure.

I have considerable affection for probabilistic arguments, and one
should resist the temptation to view these two approaches as competitors.
In fact, one purpose of this article is to explain how naturally these
perspectives can complement and reinforce one another. A coupling
may furnish the contraction estimate needed to prove a spectral gap;
a Lyapunov function may be combined with a local functional inequality
to produce a global one; and the resulting functional inequality may
then yield convergence estimates which are robust to features of the
process that would be awkward to track pathwise. The functional viewpoint
is healthily regarded as another level at which probabilistic insight
can be organised and reused.

\subsection{Questions and Scope}

Some questions which we hope to answer for the reader here include
\begin{itemize}
\item What is meant by a `functional inequality', and in what ways should
it characterise the behaviour of a Markov process?
\item In particular, what do the Poincaré and logarithmic Sobolev inequalities
say, when should one expect them to hold, and what consequences do
they yield?
\item How can one prove such an inequality using familiar probabilistic
information, including contraction and Lyapunov estimates?
\item How should a reader recognise the right inequality for a particular
convergence question and find an appropriate route into the literature?
\end{itemize}
The discussion is deliberately selective. The Poincaré inequality
(PI) and logarithmic Sobolev inequality (LSI) will be distinguished
as the main protagonists of our exposition, due to their ubiquity
and wide scope of utility. We emphasise convergence of marginal laws
and the implications for mixing times; more refined implications around
concentration of measure and geometry appear principally when they
clarify the probabilistic meaning of PI or LSI. Most proofs are suppressed
(with appropriate reference to the literature), although a few short
calculations are included when they are particularly informative.

The chosen route into the subject is through the overdamped Langevin
diffusion on Euclidean spaces, hopefully familiar to probabilists,
statisticians, machine learners, physicists, and functional analysts
alike. In studying this particular process, the influence of the target
distribution on the Markov geometry is particularly clear, various
formal calculations become particularly explicit, and intuitions from
continuous optimisation can be adopted freely. Nevertheless, the general
principles which are illustrated by this particular process will be
seen to apply to Markov processes of much broader provenance (including
those with jumps, and those which evolve in discrete time), and the
results and associated proof techniques will typically be sufficiently
abstract that their extension beyond this context is rather direct.
With this being said, we will generally attempt to indicate instances
in which the diffusion setting is essential to a specific argument.

A somewhat less innocuous restriction implied by this choice is that
we implicitly constrain our focus to the study of reversible processes,
which offers genuine simplifications in various computations, which
are not all so easily recovered in the general case. With this being
said, it should not be understood that the functional-analytic approach
advocated for herein is inapplicable outside of the reversible setting;
indeed, when adapted appropriately (through e.g. the framework of
\emph{hypocoercivity}; see e.g. \cite{villani2009hypocoercivity,dolbeault2015hypocoercivity,brigati2026hypocoercivity,fan2026sharp}),
it is even able to deliver results which are not known to be achievable
by other means. Nevertheless, we largely omit detailed discussion
of this literature.

\subsection{Notation and Priorities}

Throughout, $\pi$ denotes a probability measure on a state space
$\mathcal{X}$, and $\left(P_{t}\right)_{t\geq0}$ denotes a $\pi$-invariant
Markov semigroup. We initialise a process with $X_{0}\sim\mu_{0}$
and write $X_{t}\sim\mu_{t}=\mu_{0}P_{t}$. Our basic aim is to quantify
the approach of $\mu_{t}$ to $\pi$ as $t\to\infty$. Unless stated
otherwise, expectations and variances are taken under $\pi$; the
notation $\mathbf{E}^{x}$ instead denotes expectation under the law
of the process started from $x$.

Two conceptual priorities influence the presentation. First, motivated
by the ubiquity of high-dimensional models and scaling limits in applications,
we will pay close attention to the scaling properties of our estimates
in high dimension, with `dimension-free' bounds serving as a gold
standard of sorts. Second, while obtaining sharp constants for the
inequalities in question is certainly welcome, we will concern ourselves
principally with obtaining bounds exhibiting the correct qualitative
dependence on pertinent features of the problem, including in particular
dimension and quantitative curvature or convexity.

We parameterise our main functional inequalities by their relaxation
rates, so that larger values of these constants correspond to faster
convergence; note that this is not entirely consistent in the literature.

\subsection{Why Take a Functional Approach?}

A familiar measure of discrepancy between probability measures is
the total variation distance, given as
\begin{align*}
\mathsf{TV}\left(\mu,\nu\right) & =\sup\left\{ \left|\mu\left(A\right)-\nu\left(A\right)\right|:A\subseteq\mathcal{X}\text{ measurable}\right\} \\
 & =\frac{1}{2}\int\left|\mu-\nu\right|\left(\mathrm{d}x\right)
\end{align*}
for a suitable interpretation of $\left|\mu-\nu\right|$. This distance
has a direct interpretation in terms of coupling and hypothesis testing,
and controls the error in expectations of bounded observables. It
is therefore a natural endpoint for many mixing arguments. With this
being said, studying the convergence of Markov chains directly in
this distance can nevertheless be difficult, as outside of settings
with particularly condensed geometry (e.g. compact or finite state
spaces), it is rare for a Markov kernel to uniformly and strictly
contract the total variation distance. As such, one often seeks to
instead exhibit a contraction for some tractable quantity which upper-bounds
the total variation distance, and which is `large enough' that a contraction
is feasible.

Probabilistic approaches commonly introduce additional structure to
overcome these difficulties. For example, suppose that one identifies
a `weight' function $\varphi\geq1$ for which, when far from equilibrium,
the Markov chain under study tends to reduce $\varphi$ on average.
One can then (following e.g. \cite{hairer2011yet,del2017stochastic,del2023stability,del2025kantorovich})
introduce a `weighted total variation' distance 
\[
\mathsf{TV}_{\varphi}\left(\mu,\nu\right)=\frac{1}{2}\int\varphi\left(x\right)\left|\mu-\nu\right|\left(\mathrm{d}x\right),
\]
which accentuates dissimilarities between $\mu$ and $\nu$ in regions
of the state space for which $\varphi$ is large. From this perspective,
there is hope for the dynamics to contract $\mathsf{TV}_{\varphi}$
uniformly in this regime. Moreover, since $\mathsf{TV}_{\varphi}\geq\mathsf{TV}$,
it is clear that bounding $\mathsf{TV}_{\varphi}$ well induces a
quantitative control on $\mathsf{TV}$.

As another example, given a distance $\mathsf{d}$ on the ambient
space $\mathcal{X}$, one can consider the induced Wasserstein distance
(see e.g. \cite{villani2009optimal,chewi2025statistical})
\[
\mathcal{T}_{\mathsf{d},1}\left(\mu,\nu\right)=\inf\left\{ \mathbf{E}_{\gamma}\left[\mathsf{d}\left(X,Y\right)\right]:\gamma\in\mathsf{Couplings}\left(\mu,\nu\right)\right\} .
\]
If the dynamics of the Markov chain under study can be coupled along
iterates so that $\mathsf{d}\left(X_{t},Y_{t}\right)$ contracts exponentially
in some typical sense, then one can exhibit a uniform contraction
for $\mathcal{T}_{\mathsf{d},1}$, and again deduce convergence estimates
as a result. (Note that for generic $\mathsf{d}$, it need not be
the case that $\mathcal{T}_{\mathsf{d},1}$-convergence implies $\mathsf{TV}$-convergence).

While both approaches can be powerful, it bears mentioning that their
meaning depends on a well-chosen weight $\varphi$ or metric $\mathsf{d}$,
which imposes some burden on the analyst to make this choice well,
and introduces a degree of subjectivity in the analysis which may
not always be desirable.

An information-theoretic viewpoint avoids making either spatial choice.
In this paradigm, one declares that $\mu\approx\nu$ when it is difficult
to distinguish a draw from one law from a draw from the other. The
Neyman--Pearson principle then directs attention to the density ratio
$\rho=\mathrm{d}\mu/\mathrm{d}\nu$, which in principle contains all
information necessary for optimal classification; see e.g. Chapter
14 of \cite{polyanskiy2025information}.

Let $\Phi:\mathbf{R}_{+}\to\mathbf{R}\cup\left\{ +\infty\right\} $
be convex, with $\Phi\left(1\right)=0$. For $\mu$ absolutely continuous
with respect to $\nu$, define the density ratio divergence (often
`$f$-divergence' or `$\Phi$-divergence'; see Chapter 7 of \cite{polyanskiy2025information})
\[
\mathsf{D}_{\Phi}\left(\mu,\nu\right)=\int\nu\left(\mathrm{d}x\right)\Phi\left(\frac{\mathrm{d}\mu}{\mathrm{d}\nu}\left(x\right)\right).
\]
Standard choices of $\Phi$ recover total variation distance ($\Phi\left(r\right)=\frac{1}{2}\left|r-1\right|$),
relative entropy ($\Phi\left(r\right)=r\log r-r+1$), the chi-squared
divergence ($\Phi\left(r\right)=\left(r-1\right)^{2}$), and the squared
Hellinger distance ($\Phi\left(r\right)=\left(\sqrt{r}-1\right)^{2}$),
among others. These divergences have different strengths and impose
different requirements on the initial law, but for many standard choices,
comparison inequalities allow their decay to be converted quantitatively
into total variation convergence.

If all goes well with our Markov process of interest, then we might
hope to identify a process-dependent constant $c\in\left(0,\infty\right)$
such that whenever $\mathsf{D}_{\Phi}\left(\mu_{0},\pi\right)<\infty$,
there holds the uniform contraction
\[
\mathsf{D}_{\Phi}\left(\mu_{t},\pi\right)\leq\exp\left(-ct\right)\mathsf{D}_{\Phi}\left(\mu_{0},\pi\right).
\]
Observe that depending on the $\Phi$ in question, finiteness of the
initial divergence can be either more or less innocuous; there are
certainly cases in which singular or `cold' (i.e. far-from-equilibrium)
initialisations cause it to be infinite, and this complicates the
analysis somewhere. We will defer discussion of this issue for now,
but remark that it can be rather important when one seeks to establish
effective mixing time bounds.

This contraction estimate is plainly an attractive conclusion, but
is rather unhelpful as an assumption. A natural way to reverse-engineer
a checkable sufficient condition is to differentiate the contraction
at $t=0$: defining an `information' functional by 
\[
\mathsf{I}_{\Phi}\left(\mu_{0},\pi\right):=\left[-\partial{}_{t}\mathsf{D}_{\Phi}\left(\mu_{t},\pi\right)\right]_{t=0},
\]
we can ask whether the energy-entropy comparison 
\[
\mathsf{I}_{\Phi}\left(\mu,\pi\right)\geq c\mathsf{D}_{\Phi}\left(\mu,\pi\right)
\]
holds for every admissible $\mu$. Applying this inequality at each
point along the semigroup and invoking Grönwall's lemma would then
give the desired exponential decay. Conversely, differentiating a
uniform contraction estimate recovers the infinitesimal inequality
under the usual regularity assumptions.

This ``reverse-Grönwall'' manoeuvre serves as a concise mnemonic for
the main inequalities which we will consider herein. Working with
inequalities of this form takes our initial dynamic question about
a Markov process and recasts it as a static comparison between a `global'
entropy functional and its `local' instantaneous dissipation. 

\section{The Langevin Diffusion and its Functional Inequalities}

In the context of the overdamped Langevin diffusion, the relationship
between these entropic and information functionals becomes particularly
explicit, and readily gives rise to the primary functional inequalities
which we seek to study.

\subsection{The Overdamped Langevin Diffusion}

\subsubsection{Basics}

Let $\pi$ be a smooth, strictly positive probability density on $\mathbf{R}^{d}$
and consider the overdamped Langevin diffusion
\[
\mathrm{d}X_{t}=\nabla\log\pi\left(X_{t}\right)\,\mathrm{d}t+\sqrt{2}\,\mathrm{d}W_{t}.
\]
Good general references on this process include \cite{bakry2014analysis,pavliotis2014stochastic,van2016probability}.
Under standard conditions this process is non-explosive, ergodic,
and reversible with invariant law $\pi$; we suppress the analytic
assumptions needed to justify these statements and the calculations
below, and so the reader may initially regard all differentiations
and integrations by parts as formal, bearing in mind that they can
generally be justified granted sufficient smoothness and sufficiently
moderate growth of all functions and densities in question. Define
the semigroup $\left(P_{t}\right)_{t\geq0}$ by $P_{t}F\left(x\right)=\mathbf{E}^{x}\left[F\left(X_{t}\right)\right]$. 

While this semigroup nominally contains the distributional information
of interest, it is rarely available in closed form. As such, we often
work instead with its (considerably more tractable) infinitesimal
description. The infinitesimal generator $\mathcal{L}$ of the semigroup
is defined pointwise by
\[
\mathcal{L}F\left(x\right)=\lim_{t\to0^{+}}\frac{P_{t}F\left(x\right)-F\left(x\right)}{t}
\]
for all $F$ such that the definition makes sense. Formally, one may
write $P_{t}$ as $\exp\left(t\mathcal{L}\right)$ and express its
evolution via $\partial_{t}P_{t}=\mathcal{L}P_{t}=P_{t}\mathcal{L}$.
Itô's formula gives the explicit expression
\[
\mathcal{L}F\left(x\right)=\left\langle \nabla\log\pi\left(x\right),\nabla F\left(x\right)\right\rangle +\Delta F\left(x\right).
\]
The associated Dirichlet energy form is the bilinear form $\mathcal{E}$,
defined by
\begin{align*}
\mathcal{E}\left(F,G\right) & =\mathbf{E}\left[F\left(-\mathcal{L}G\right)\right]=\mathbf{E}\left[\left\langle \nabla F,\nabla G\right\rangle \right],
\end{align*}
where both expectations are computed under $\pi$, and where the latter
equality follows for suitable functions $F$, $G$ by integration
by parts. The symmetry of this form is a functional-analytic reflection
of the dynamical reversibility property (`detailed balance') of the
underlying dynamics. We will see that this form provides an analytic
language for computing and comparing the dissipation of different
pertinent functionals.

\subsubsection{Entropy Dissipation Computations}

Return now to our `reverse-Grönwall' programme: we seek to make statements
of the form 
\[
\mathsf{D}_{\Phi}\left(\mu_{t},\pi\right)\lesssim_{t}\mathsf{D}_{\Phi}\left(\mu_{0},\pi\right)
\]
which hold for generic $\mu_{0}$. As a first step, reversibility
of the semigroup allows one to compute that the density of $\mu_{t}$
with respect to the target $\pi$ evolves as
\[
\frac{\mathrm{d}\mu_{t}}{\mathrm{d}\pi}=P_{t}\frac{\mathrm{d}\mu_{0}}{\mathrm{d}\pi},
\]
and so upon defining the `$\Phi$-entropy' functional $\mathsf{Ent}_{\Phi}$
by
\[
\mathsf{Ent}_{\Phi}\left(F\right)=\mathbf{E}\left[\Phi\circ F\right]-\Phi\left(\mathbf{E}\left[F\right]\right)\geq0,
\]
one sees that
\[
\mathsf{D}_{\Phi}\left(\mu_{t},\pi\right)=\mathsf{Ent}_{\Phi}\left(P_{t}\rho_{0}\right),\qquad\rho_{0}=\frac{\mathrm{d}\mu_{0}}{\mathrm{d}\pi}.
\]
Note that $\mathbf{E}\left[P_{t}\rho_{0}\right]=1$ for all $t$,
and hence $\Phi\left(\mathbf{E}\left[P_{t}\rho_{0}\right]\right)=0$
for all time; this simplifies some of the computations which follow.

From this perspective, we see that it is sufficient to establish that
for a sufficiently large class of functions $F$, one can control
\[
\mathsf{Ent}_{\Phi}\left(P_{t}F\right)\lesssim_{t}\mathsf{Ent}_{\Phi}\left(F\right).
\]
In making this passage, we transition from reasoning about probability
measures to reasoning directly about functions, which can often simplify
matters.

Given this equivalence, we now compute time derivatives as anticipated,
calculating that (assuming sufficient smoothness for $\Phi$ and $F$)
\begin{align*}
-\partial_{t}\mathsf{Ent}_{\Phi}\left(P_{t}F\right) & =-\partial_{t}\left(\mathbf{E}\left[\Phi\left(P_{t}F\right)\right]-\Phi\left(\mathbf{E}\left[P_{t}F\right]\right)\right)\\
 & =\mathbf{E}\left[\Phi^{\prime}\left(P_{t}F\right)\left(-\mathcal{L}P_{t}F\right)\right]\\
 & =\mathcal{E}\left(\Phi^{\prime}\left(P_{t}F\right),P_{t}F\right)\\
 & =\mathbf{E}\left[\Phi^{\prime\prime}\left(P_{t}F\right)\left\Vert \nabla\left(P_{t}F\right)\right\Vert ^{2}\right],
\end{align*}
which is often described as a sort of `$\Phi$-Information' functional.

On the basis of these calculations, one might be emboldened to say
that the probability measure $\pi$ satisfies a `$\Phi$-Sobolev inequality
with constant $c\geq0$' if for all suitable $G$, there holds the
estimate
\[
\mathbf{E}\left[\Phi^{\prime\prime}\left(G\right)\left\Vert \nabla G\right\Vert ^{2}\right]\geq c\mathsf{Ent}_{\Phi}\left(G\right).
\]
If this holds, then we compute that
\begin{align*}
-\partial_{t}\mathsf{Ent}_{\Phi}\left(P_{t}F\right) & =\mathbf{E}\left[\Phi^{\prime\prime}\left(P_{t}F\right)\left\Vert \nabla\left(P_{t}F\right)\right\Vert ^{2}\right]\geq c\mathsf{Ent}_{\Phi}\left(P_{t}F\right),
\end{align*}
whence Grönwall's lemma yields that
\[
\mathsf{Ent}_{\Phi}\left(P_{t}F\right)\leq\exp\left(-ct\right)\mathsf{Ent}_{\Phi}\left(F\right)
\]
for all $F$, and hence that 
\[
\mathsf{D}_{\Phi}\left(\mu_{t},\pi\right)\leq\exp\left(-ct\right)\mathsf{D}_{\Phi}\left(\mu_{0},\pi\right)
\]
for all suitable $\mu_{0}$, as desired.

We have therefore converted exponential contraction of the Markov
dynamics into a static inequality of Sobolev type, comparing the global
variability of $F$ under $\pi$ to its local variability. The difficult
part is now to determine when such an inequality holds and how the
associated constant $c$ depends on the problem. Chafaï's survey \cite{chafai2004entropies}
develops this calculation for a broad family of entropies generated
by convex functions $\Phi$. 

\subsection{Poincaré and Log-Sobolev: The Main Protagonists}

Here, we narrow our attention to the two cases which recur most often
in the study of Markov processes.

\subsubsection{Definitions and Dynamical Meaning}

We start by writing down the form of these `$\Phi$-Sobolev' inequalities
for some important choices of $\Phi$.
\begin{itemize}
\item With $\Phi\left(r\right)=\left(r-1\right)^{2}$, some simplifications
yield the `Poincaré inequality' (PI) 
\[
\mathbf{E}\left[\left\Vert \nabla G\right\Vert ^{2}\right]\geq\gamma\operatorname{Var}_{\pi}\left(G\right),
\]
i.e. the variance of a function is controlled by its expected squared
gradient. Equivalently, $\gamma$ is a lower bound on the non-zero
spectrum of the linear operator $-\mathcal{L}$; the largest admissible
value is known as the `spectral gap' of the process. In dynamical
terms, it is the $\mathrm{L}^{2}\left(\pi\right)$ relaxation rate
of the process when the optimal constant is used. Because the chi-squared
divergence is a \emph{squared} $L^{2}$ distance between densities,
it thereby contracts at the rate $2\gamma$.
\item With $\Phi\left(r\right)=r\log r-r+1$, one instead finds the so-called
`logarithmic Sobolev inequality' (LSI)
\[
\mathbf{E}\left[\frac{\left\Vert \nabla G\right\Vert ^{2}}{G}\right]\geq2\lambda\mathsf{Ent}\left(G\right),
\]
where $\mathsf{Ent}\left(G\right)=\mathbf{E}\left[G\log G\right]-\mathbf{E}\left[G\right]\log\mathbf{E}\left[G\right]$
for nonnegative $G$. This is often transformed into the equivalent
forms (replacing $G\leftarrow G^{2}$, $G\leftarrow e^{G}$ respectively)
\[
\mathbf{E}\left[\left\Vert \nabla G\right\Vert ^{2}\right]\geq\frac{\lambda}{2}\mathsf{Ent}\left(G^{2}\right),\qquad\mathbf{E}\left[\left\Vert \nabla G\right\Vert ^{2}e^{G}\right]\geq2\lambda\mathsf{Ent}\left(e^{G}\right).
\]
Correspondingly, the relative entropy $\mathsf{KL}\left(\mu_{t},\pi\right)$,
also called the `Kullback--Leibler divergence', contracts at rate
$2\lambda$.
\end{itemize}
The LSI with constant $\lambda$ implies a PI with $\gamma\geq\lambda$.
Morally, this is a consequence of the fact that as $r\to1$, there
is the equivalence $r\log r-r+1\asymp\frac{1}{2}\left(r-1\right)^{2}$;
rigorously, this can be proven by considering $G_{\varepsilon}=1+\varepsilon G$
in the LSI and sending $\varepsilon\to0^{+}$. On the other hand,
we will see that in general, the LSI is strictly stronger than the
PI; that is, there are measures for which PI holds with $\gamma>0$
while LSI fails for every $\lambda>0$.

One might then reasonably ask why these two inequalities are such
key figures in the literature. At an informal level, they delimit
a particularly useful range of behaviours for high-dimensional models:
\begin{enumerate}
\item Weakenings of PI generally correspond to Markov processes with slower-than-exponential
relaxation to equilibrium. As such, one loses the ability to quantify
the convergence properties of the process with a single exponential
rate, and some finer description of the convergence profile becomes
necessary.
\item Strengthenings of LSI generally entail `better-than-Gaussian' behaviour,
in terms of tail behaviour and convergence behaviour when far from
equilibrium. While such behaviours are often desirable for any fixed
problem, they stand as an obstruction to dimension-free behaviour
and tensorisation properties, essentially because, as a consequence
of the Central Limit Theorem, `better-than-Gaussian' behaviour can
scarcely continue indefinitely for product measures.
\end{enumerate}
This is not to say that PI and LSI are the only worthwhile inequalities,
but experience certainly dictates that they strike an especially favourable
balance vis-à-vis broad applicability and user-friendliness.

\subsubsection{First Examples}

In the context of high-dimensional probability, one often benefits
from first considering Gaussian measures. Indeed, letting $\pi=\mathcal{N}\left(0_{d},\mathsf{C}\right)$
denote a Gaussian measure and taking $\rho=\lambda_{\mathrm{min}}\left(\mathsf{C}^{-1}\right)$,
one can show that $\pi$ satisfies the logarithmic Sobolev inequality
with constant $\geq\rho$ (and hence PI holds with the same constant).

While the specifics of the Gaussian measure are convenient for explicit
calculations, it happens that this particular result is actually rather
robust: assuming $\pi$ to be $\rho$-strongly log-concave in the
sense that uniformly in $x$, one has the matrix ordering
\[
\nabla^{2}\left(-\log\pi\right)\left(x\right)\succeq\rho\mathbf{I}_{d},
\]
then the same guarantee persists, i.e.
\[
\mathbf{E}\left[\left\Vert \nabla G\right\Vert ^{2}\right]\geq\frac{\rho}{2}\mathsf{Ent}\left(G^{2}\right).
\]
There are numerous ways to establish this result, several of which
will be discussed later on as applications of broader principles.
Among other interpretations, this is the basic curvature criterion
of the Bakry--Émery `$\Gamma$-calculus', which can be used to prove
this LSI and much more. The interested reader is referred to the monograph
\cite{bakry2014analysis} of Bakry, Gentil, Ledoux on `Analysis and
Geometry of Markov Diffusion Operators', which can be seen as attempting
to pursue `all possible consequences' of this particular assumption.
The conclusions are certainly powerful, but one quickly realises that
uniform positive curvature is a rather rigid assumption. Part of the
appeal of the methods which will be treated in the sequel is that
they can continue to supply useful estimates even after this pointwise
condition has failed.

For another grounding example, take $\pi$ to be the unit-rate exponential
distribution on $\mathbf{R}^{+}$, i.e. $\pi\left(x\right)=\exp\left(-x\right)\mathbf{1}_{x>0}$.
For this distribution, the Poincaré inequality holds with constant
$1/4$, i.e. it holds that
\[
4\mathbf{E}\left[\left\Vert \nabla G\right\Vert ^{2}\right]\geq\operatorname{Var}_{\pi}\left(G\right).
\]
One approach to establishing this involves studying the associated
eigenproblem for the generator $\mathcal{L}$ more-or-less directly.
On the dynamical side this corresponds to the Langevin diffusion on
$\mathbf{R}_{+}$ with prevailing dynamic
\[
\mathrm{d}X_{t}=-\mathrm{d}t+\sqrt{2}\,\mathrm{d}W_{t}
\]
and understood to have a reflecting boundary at the origin. Considering
test functions of the form $G\left(x\right)=\exp\left(\lambda x\right)$
allows one to establish sharpness of the constant $\gamma=1/4$. On
the other hand, the same collection of test functions can be used
to exhibit that the LSI fails for $\pi$ for every positive constant.
We thus see quickly and cleanly that the two inequalities are genuinely
separated; intuition for this separation will hopefully become yet
clearer over the course of this article.

Thus far, we have a small collection of highly-structured exemplar
probability measures which we know to satisfy the inequalities of
interest. Our next goal is to demonstrate that these examples can
easily be broadcast into a much richer and more diverse collection
of examples, using a sort of `algebra' for functional inequalities.

\subsection{Transfer Principles for Functional Inequalities}

How can we prove that such inequalities can possibly hold for more
interesting probability measures $\pi$? It will not be viable to
solve case after case by explicit computations, and so we are driven
to understand how we can combine existing probability measures in
a `Sobolev-compatible' manner. In particular, tensorisation, Lipschitz
transport, convolution, bounded changes of measure, and sufficiently
overlapping mixtures will be seen to transform the constants in simple
and interpretable ways.

In what follows, I will focus my initial discussion on the Poincaré
inequality, but most developments will also be valid for the logarithmic
Sobolev inequality, and occasionally, even for more general $\Phi$-Sobolev
inequalities; I will add parenthetical remarks to this effect where
appropriate.

\subsubsection{Tensorisation}

Suppose that for $1\leq i\leq N$, $\mu_{i}$ is a probability measure
which satisfies a Poincaré inequality with constant $\gamma_{i}>0$,
and define $\mu=\bigotimes_{1\leq i\leq N}\mu_{i}$. It then holds
that $\mu$ satisfies a Poincaré inequality with constant $\gamma=\min\left\{ \gamma_{i}:1\leq i\leq N\right\} $.
The same conclusion holds verbatim for the LSI, and for a structured
sub-class of other $\Phi$-Sobolev inequalities; see the paper \cite{chafai2004entropies}
of Chafaï for further details.

Informally, this conveys that from the perspective of functional inequalities,
a product measure in this class (or equivalently, the rate of mixing
of a product-form Markov process) is `no worse than its worst one-dimensional
component'. On the basis of this result, one can easily exhibit high-dimensional
measures which enjoy the Poincaré inequality simply by taking products. 

\subsubsection{Transport}

Suppose that $\mu$ is a probability measure on $\mathbf{R}^{d}$
which satisfies a Poincaré inequality with constant $\gamma>0$, and
let $T:\mathbf{R}^{d}\to\mathbf{R}^{d^{\prime}}$ be an $L$-Lipschitz
mapping for some $L<\infty$. Defining then the pushforward measure
$\pi:=T_{\#}\mu$ (i.e. $\pi\left(A\right)=\mu\left(\left\{ x\in\mathbf{R}^{d}:T\left(x\right)\in A\right\} \right)$),
it holds that $\pi$ satisfies a Poincaré inequality with constant
at least $\gamma L^{-2}$. Again, the same conclusion holds verbatim
for the LSI, and indeed, for various other inequalities involving
gradient norms in their upper bounds.

Simple applications of this principle include that these inequalities
are preserved exactly under orthogonal changes of basis, preserved
up to constants under general linear changes of basis, and can be
improved upon marginalisation (where e.g. $T:\left(x,y\right)\mapsto x$).
More advanced applications involve studying more varied nonlinear
mappings for which Lipschitz regularity can be established.

\begin{figure}[h]
\caption{A Lipschitz transport can transfer a reference measure to a more geometrically
complicated target, while controlling pairwise distances.}

\centering{}\includegraphics[scale=0.2]{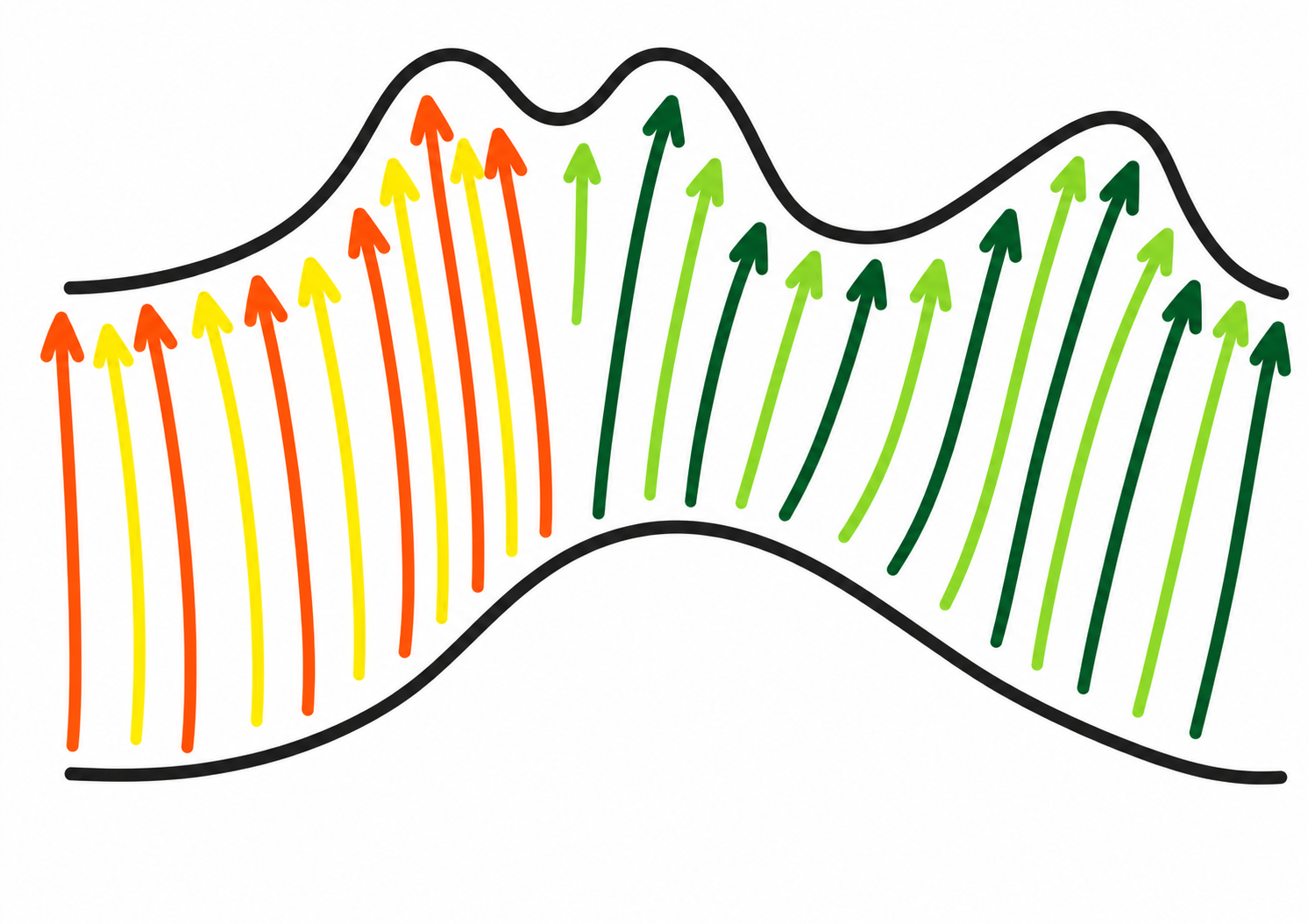}
\end{figure}

Indeed, this principle enables the general possibility of exhibiting
functional inequalities for measures which are `metric contractions'
of measures which are already known to be well-behaved. This is particularly
straightforward to apply in one dimension, where the monotone transport
map between two regular probability measures is uniquely identified
in terms of cumulative distribution functions. In the multivariate
case, one has to work slightly harder in exhibiting specific maps
which are known to have good Lipschitz regularity; contemporary constructions
of such maps are based on deep ideas from optimal transport \cite{caffarelli2000monotonicity,kolesnikov2013sobolev,fathi2020proof},
stochastic dynamics and the heat flow \cite{kim2012generalization,neeman2022lipschitz,mikulincer2024brownian,conforti2025coupling},
and beyond. Based on these constructions, one can exhibit quantitative
functional inequalities (and more) for strongly log-concave measures,
Gaussian mixtures with compact mixing measures, smooth and bounded
perturbations of Gaussian measures, and a small handful of other cases.
Constructive Lipschitz transport is now an active subject in its own
right, bolstered in part by a timely connection to thrilling contemporary
developments in generative modelling by measure transport; see \cite{song2021score,lipman2023flow,albergo2025stochastic}
for some representative references from the machine learning literature.

By contrast with some of the other tools outlined in this section,
extending this approach to the study of non-diffusive Markov processes
seems to be rather more difficult. One fascinating discrete example
is the Poisson transport map of López-Rivera and Shenfeld \cite{lopez2025poisson},
which is `sufficiently Lipschitz' in character to be able to transfer
functional inequalities from Poisson point processes to ultra-log-concave
measures on the natural numbers. Identifying further examples in more
varied settings would be of substantial interest.

\subsubsection{Convolution}

On the basis of these elementary results, one can already prove a
simple but non-trivial result about convolution measures. For $1\leq i\leq N$,
let $X_{i}\sim\mu_{i}$ independently, where $\mu_{i}$ is a probability
measure on $\mathbf{R}^{d}$ which satisfies a Poincaré inequality
with constant $\gamma_{i}>0$, and define $Y=\sum_{1\leq i\leq N}X_{i}$.
The joint law of $X=\left\{ X_{i}:1\leq i\leq N\right\} $ is of product
form and the summation map $X\mapsto Y$ is $\sqrt{N}$-Lipschitz,
giving the immediate estimate $\gamma_{Y}\geq N^{-1}\min\left\{ \gamma_{i}:1\leq i\leq N\right\} $.
A sharpened argument of Chafaï (which is still relatively elementary;
see Section 3.2 of \cite{chafai2004entropies}) improves this to $\gamma_{Y}\geq\left(\sum_{1\leq i\leq N}\gamma^{-1}_{i}\right)^{-1}$,
which is seen to be considerably more accurate in the case of inhomogeneous
summands. Again, the same conclusions hold verbatim for the LSI and
a structured sub-class of other $\Phi$-Sobolev inequalities. As one
might imagine, these results can be useful for quantitative approaches
to and reinforcements of the Central Limit Theorem, wherein convolution
and `Gaussian-like' behaviour (generously interpreted) play a central
role.

\subsubsection{Bounded Change-of-Measure}

The preceding examples have all involved transformations which are,
informally speaking, regular `at the level of samples'. The following
example involves a transformation which is instead regular `at the
level of densities'. More precisely, suppose that $\mu$ is a probability
measure on $\mathbf{R}^{d}$ which satisfies a Poincaré inequality
with constant $\gamma>0$, and let $\pi$ be a probability measure
on $\mathbf{R}^{d}$ which is `strongly mutually absolutely continuous'
with respect to $\mu$, in the sense that $\sup\frac{\mathrm{d}\pi}{\mathrm{d}\mu}$
and $\sup\frac{\mathrm{d}\mu}{\mathrm{d}\pi}$ are both finite. Writing
\[
\kappa:=\left(\sup\frac{\mathrm{d}\pi}{\mathrm{d}\mu}\right)\left(\sup\frac{\mathrm{d}\mu}{\mathrm{d}\pi}\right)\geq1
\]
 for the `condition number' of the pair of measures, it then follows
that $\pi$ satisfies a Poincaré inequality with constant $\geq\gamma\kappa^{-1}$.
The same conclusion holds verbatim for the LSI, all $\Phi$-Sobolev
inequalities, and even for many energy-entropy inequalities corresponding
to very different geometries.

\begin{figure}[h]
\caption{Two views of a bounded perturbation: i) a non-strongly convex potential
might be close to some strongly convex potential in $\sup$-norm,
or ii) a density of interest might be sandwiched well between two
multiples of a well-behaved reference density.}

\centering{}\includegraphics[scale=0.14]{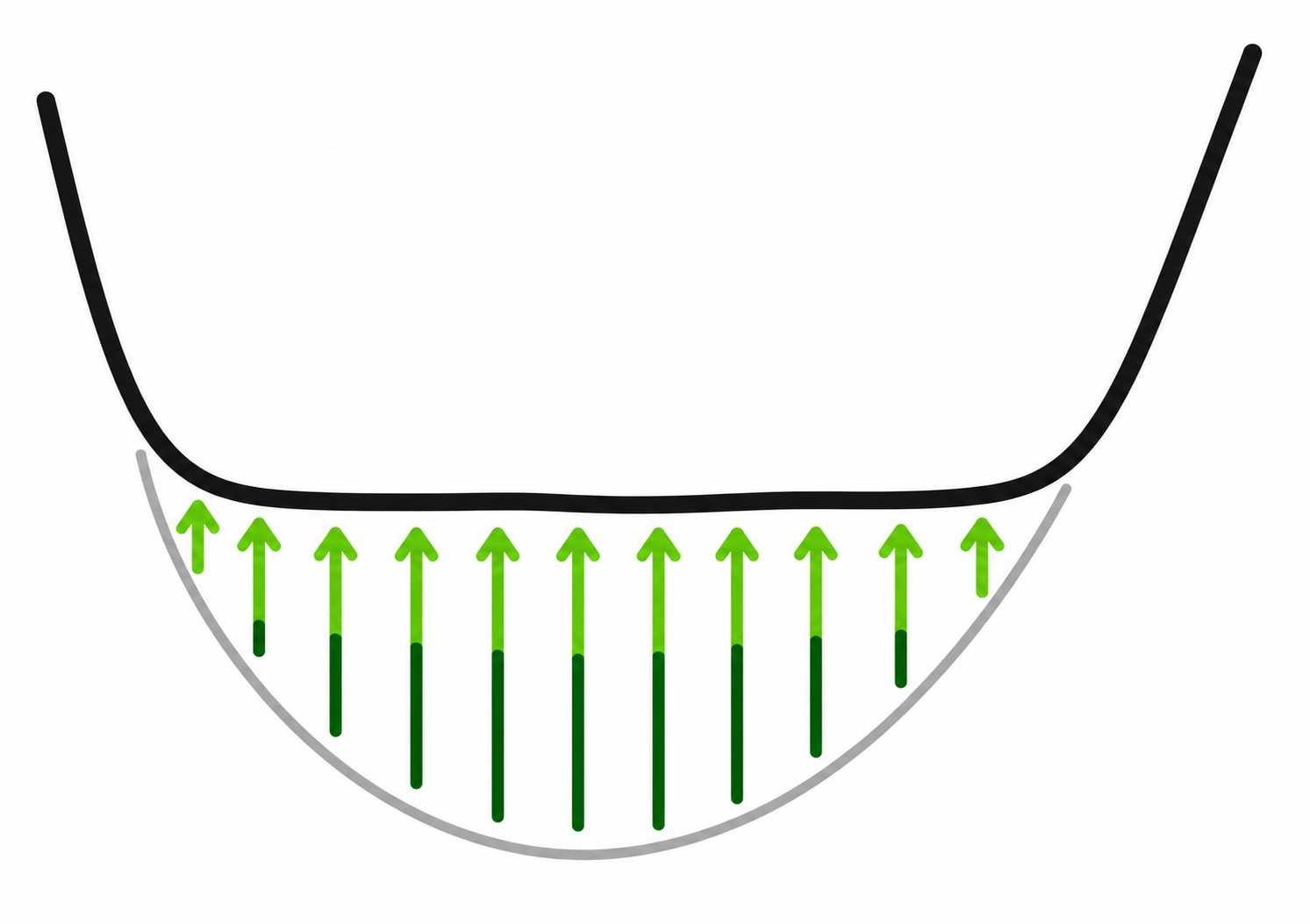}\hfill{}\includegraphics[scale=0.14]{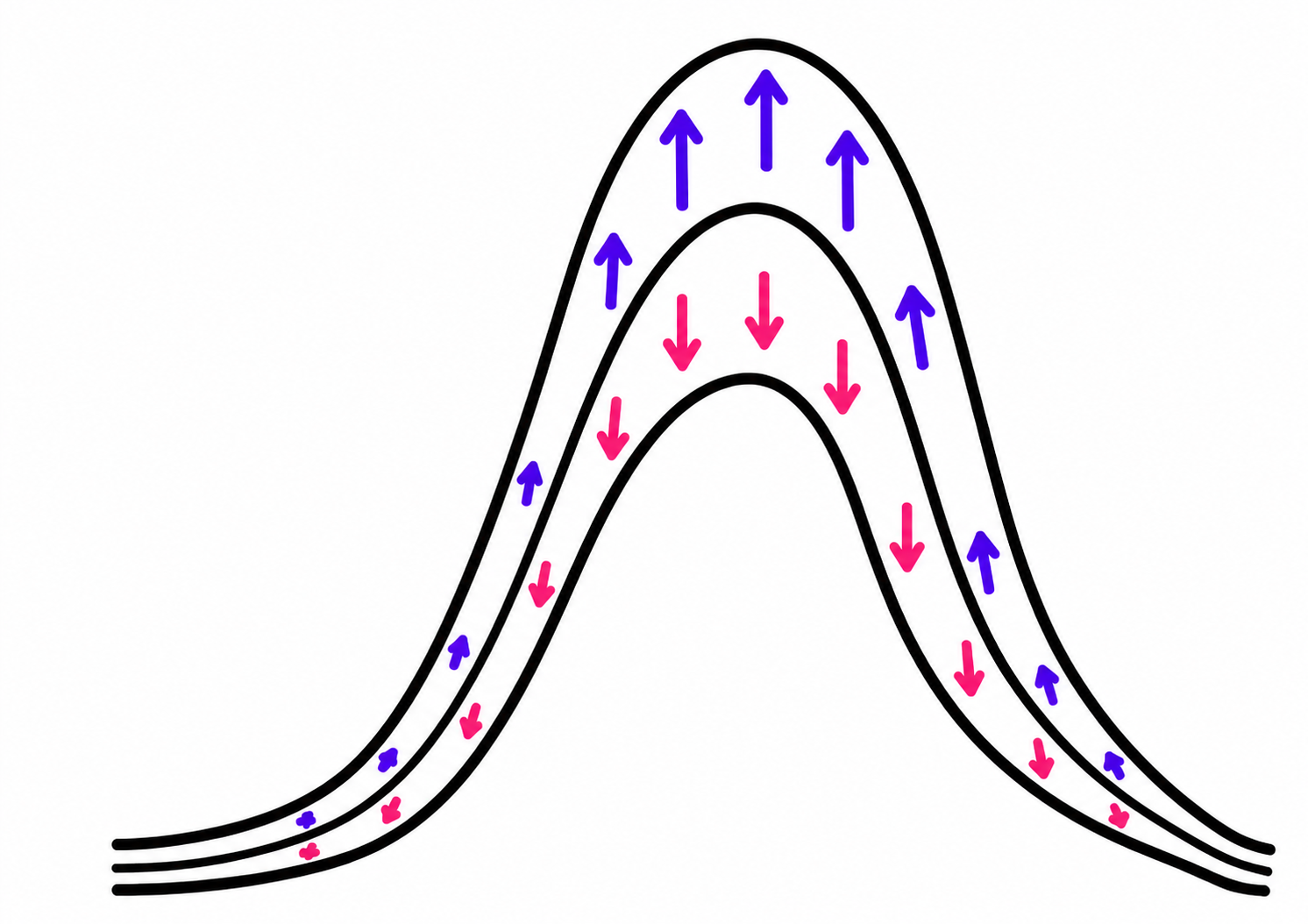}
\end{figure}

This shows that even probability measures with rather rough densities
can still enjoy good functional inequalities, provided that the roughness
is well-confined in a uniform sense. Moreover, it highlights that
while convexity properties are favourable for establishing these inequalities,
they are really quite far from necessary. 

As a simple example, for `double-well' potentials of the form $U_{\lambda}\left(x\right)=x^{4}-\lambda x^{2}$
with $\lambda>0$, one can write
\[
U_{\lambda}\left(x\right)=\lambda x^{2}+\max\left\{ 0,x^{4}-2\lambda x^{2}\right\} +\min\left\{ 0,x^{4}-2\lambda x^{2}\right\} ,
\]
observe that the first two terms are collectively strongly convex
with curvature $2\lambda$, and that the third term is bounded between
$-\lambda^{2}$ and $0$. One then observes that $\pi\propto\exp\left(-U_{\lambda}\right)$
enjoys both the LSI and PI with constant at least $2\lambda\exp\left(-\lambda^{2}\right)$.
More broadly, the same flavour of argument allows one to show that
many ``asymptotically log-concave'' measures will enjoy the PI with
some positive constant; see \cite{chen1997estimates} for a particularly
general result of this flavour.

This estimate and its extensions are typically referred to as the
(Deuschel-)Holley-Stroock principle, following early contributions
\cite{holley1987logarithmic,deuschel1990hypercontractivity}. Anecdotally,
this approach has something of a reputation for being rather limited,
since bounded perturbations are a priori also apparently quite limited.
Nevertheless, consider that if one is interested in establishing functional
inequalities for a specific $\pi$, one has flexibility in choosing
which `nice' reference $\mu$ to consider, and trading off proximity
to $\pi$ against the strength of the functional inequality enjoyed
by $\mu$. Taking this perspective, the principle can be rather more
powerful.

\subsubsection{Mixtures}

A complementary question is whether mixtures of well-behaved distributions
remain well behaved. Observe that some assumptions are clearly required
here: from the perspective of functional inequalities, Gaussian measures
are about as nice as one can hope for, but through the lens of functional
inequalities, well-separated Gaussian mixtures can be quite awful.
In this regard, it is not surprising that some sort of `overlap' condition
enables positive consequences.

More precisely, a result of Chen--Chewi--Niles-Weed \cite{chen2021dimension}
establishes that if $\mathcal{P}=\left\{ P_{z}:z\in\mathcal{Z}\right\} $
is a family of probability measures on $\mathcal{X}=\mathbf{R}^{d}$
which each satisfy a Poincaré inequality with constant $\geq\gamma$,
and if the family $\mathcal{P}$ is of `finite diameter in the $\chi^{2}$
sense', i.e. with $\chi^{2}$ denoting the chi-squared divergence,
there exists a finite $\bar{\chi^{2}}<\infty$ such that
\[
\sup\left\{ \chi^{2}\left(P_{z},P_{z^{\prime}}\right):z,z^{\prime}\in\mathcal{Z}\right\} \leq\bar{\chi^{2}},
\]
then for any probability measure $\mu$ over $\mathcal{Z}$, the mixture
$\pi=\mu P$ also satisfies a Poincaré inequality with explicit constant
depending only on $\gamma$ and $\bar{\chi^{2}}$. Their work also
includes results under weaker moment-type conditions on the pairwise
chi-squared divergences, similar results for LSI, and more besides.
In each case, the mixtures $\pi$ can be far from log-concave, but
still exhibit good connectivity properties.

This result and its proof reflect cleanly that while the component
inequalities control variation within each component, the divergence
bound is needed in order to control the fluctuations of expectations
across components; one sees that well-separated mixtures fail (at
least quantitatively) precisely at the second stage.

\subsection{Probabilistic Routes to Functional Inequalities}

While the preceding transfer principles are powerful and useful, they
ultimately rely on some initial functional inequality being proven
in the first place. We now highlight two routes to doing so which
are particularly natural to readers familiar with probabilistic ergodic
theory: contraction by coupling, and Lyapunov drift conditions.

\begin{figure}[h]
\caption{Two probabilistic routes to functional inequalities: coupling Markov
dynamics to bring the two copies closer together in space (left),
and using a drift condition to drive arbitrary initialisations towards
a sublevel set on which the process behaves simply (right).}

\centering{}\includegraphics[scale=0.14]{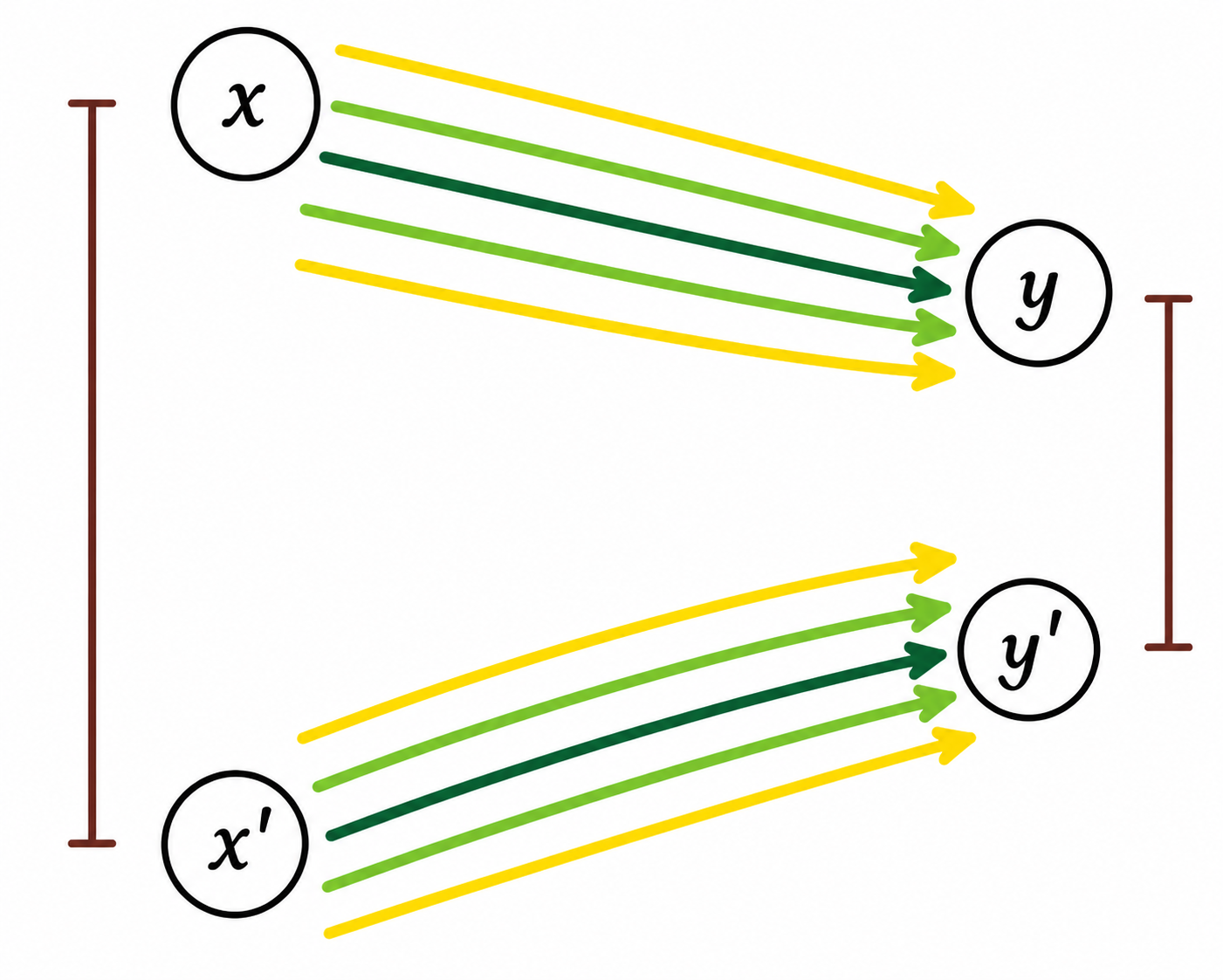}\hfill{}\includegraphics[scale=0.14]{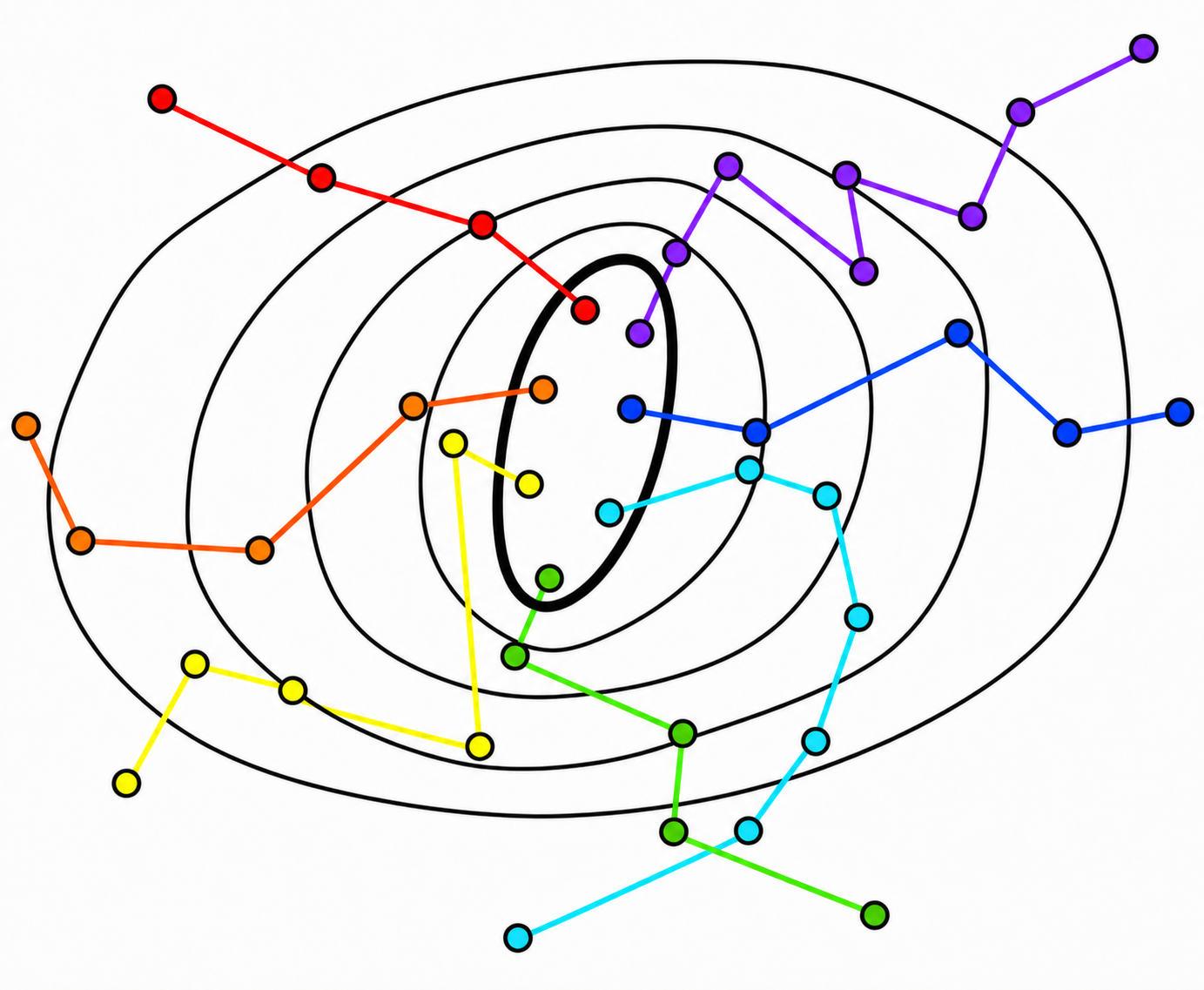}
\end{figure}

\subsubsection{Coupling and Curvature}

A useful principle, appearing in work of Mu-Fa Chen \cite{chen2004markov,chen2005eigenvalues},
Feng-Yu Wang \cite{chen1994application}, and Yann Ollivier \cite{ollivier2009ricci},
among others, is that a contraction by coupling (see e.g. \cite{lindvall2002lectures,den2012probability}
for overviews) can imply a positive spectral gap, i.e. a Poincaré
inequality. Under suitable regularity and compatibility assumptions,
suppose that a reversible (in fact, normality suffices) Markov semigroup
$\left\{ P_{t}:t\geq0\right\} $ enjoys a Wasserstein contraction
condition of the form
\[
\forall x,x^{\prime},\qquad\mathcal{T}_{\mathsf{d},1}\left(P_{t}\left(x,\cdot\right),P_{t}\left(x^{\prime},\cdot\right)\right)\leq\exp\left(-\rho t\right)\mathsf{d}\left(x,x^{\prime}\right)
\]
for some $\rho>0$. It then follows that the same semigroup enjoys
a Poincaré inequality with constant $\gamma\geq\rho$. 

While the easiest applications of this principle involve taking $\mathsf{d}$
to be e.g. the standard Euclidean metric, more advanced applications
incorporate the idea that the choice of metric $\mathsf{d}$ is part
of the argument rather than part of the model, leaving room to improve
the estimate by designing $\mathsf{d}$ carefully. In the context
of the Langevin diffusion, this degree of freedom is of particular
utility for non-uniformly convex potentials, where careful design
of both coupling and metric can deliver estimates which would not
be accessible by more naive methods. A typical example in this literature
involves potentials satisfying curvature bounds of the form
\[
\nabla^{2}U\left(x\right)\succeq\begin{cases}
-L\mathbf{I}_{d} & \left\Vert x\right\Vert <r\\
M\mathbf{I}_{d} & \left\Vert x\right\Vert \geq r,
\end{cases}
\]
and working with metrics of the form $\mathsf{d}\left(x,y\right)=f\left(\left\Vert x-y\right\Vert \right)$
for a carefully-designed concave function $f$; see \cite{eberle2016reflection}
for a detailed development of this approach.

While Wasserstein contraction implies an estimate on the spectral
gap, in general, it does not imply the corresponding entropy contraction,
or modified logarithmic Sobolev inequality, without additional hypotheses,
as confirmed by a recent work of Münch \cite{munch2023ollivier}.
Recent developments on the so-called `Peres--Tetali conjecture'\footnote{Anecdotally, it seems to be the case that while Peres and Tetali might
have each been interested in the potential implications from metric
contraction to entropy contraction, they did not particularly reach
consensus on whether this implication ought to follow without additional
hypotheses; see Tetali's recollections in this video (\href{https://youtu.be/J2l8s4roWis?si=3r1t7rAODkSqGCfQ}{https://youtu.be/J2l8s4roWis?si=3r1t7rAODkSqGCfQ})
at around 40 minutes.} show that a non-negative sectional-curvature condition, interpretable
as pathwise non-expansion for suitable backward dynamics, can bridge
this gap. Roughly speaking, if one has contraction-on-average for
the forward dynamics, complemented by an almost-sure non-expansion
for the backward dynamics, then the work of Caputo--Münch--Salez
\cite{caputo2025entropy} proves that an entropy contraction ensues,
and supplies a quantitative estimate for the associated contraction
rate. Beyond the plainly useful practical consequences for estimating
convergence rates, this example supplies some useful intuition about
the qualitative separation between variance contraction and entropy
contraction, whereby the latter is typically demanding of more regular,
uniform convergence behaviour. This remarkable result also admits
a convenient extension \cite{caputo2026entropy} to `multi-block'
systems, as encountered in the context of Gibbs sampling / `Glauber
dynamics'.

A related approach involves looking at `typical' rather than `uniform'
curvature. Returning to the specific context of the Langevin diffusion,
assume that for all $x$, one can bound
\[
\nabla^{2}\left(-\log\pi\right)\left(x\right)\succeq\rho\left(x\right)\mathbf{I}_{d}
\]
with $\rho\left(x\right)>0$ globally; this can be understood as a
mild reinforcement of strict (but not strong) log-concavity. A result
of Veysseire \cite{veysseire2010harmonic,veysseire2011improved} then
shows that $\pi$ enjoys a Poincaré inequality with estimate of `harmonic
mean' type, i.e. 
\[
\gamma\geq\mathbf{E}\left[\rho\left(X\right)^{-1}\right]^{-1},
\]
provided that the expectation in question is finite. Later work situates
the proof of this result alongside related methods of `intertwining'
and `generalised Brascamp-Lieb' inequalities; see e.g. \cite{arnaudon2018intertwinings,bonnefont2021intertwinings,huguet2022intertwining,bonnefont2024note}
for an entry to these techniques.

While a sharper `arithmetic mean' estimate of the form $\gamma\geq\mathbf{E}\left[\rho\left(X\right)\right]$
would be appealing, it is known that it cannot hold in general, even
for strictly convex potentials. Recent work of Cattiaux-Fathi-Guillin
\cite{cattiaux2022self} obtains related bounds involving some additional
assumptions on the regularity of the local curvature $\rho$.

\subsubsection{Lyapunov Conditions}

It is well-known that the `Drift and Minorisation' approach to stochastic
stability provides tools for establishing exponential ergodicity;
see the authoritative textbook of \cite{meyn2009markov} for full
developments. For a more localised introduction, assume that 
\begin{enumerate}
\item for some $V:\mathbf{R}^{d}\to\left[1,\infty\right)$, there holds
for some $\lambda>0$, $K>0$ the drift (or `Lyapunov') condition
\[
\mathcal{L}V\leq-\lambda V+K,
\]
and that 
\item for some `inflation' factor $\iota>1$, there holds a `minorisation'
on $C_{\iota}:=\left\{ x\in\mathbf{R}^{d}:V\left(x\right)\leq\iota K/\lambda\right\} $,
i.e. there exists a $\tau>0$, $\varepsilon\in\left(0,1\right)$,
and a probability measure $\alpha$ such that
\[
x\in C_{\iota}\implies P_{\tau}\left(x,\cdot\right)\geq\varepsilon\alpha\left(\cdot\right)
\]
as measures.
\end{enumerate}
Informally, the first condition suggests that the process will rapidly
enter (or return) into the set $C_{\iota}$, wherein the second condition
implies that upon entering into this set, different copies of the
process can be made to coalesce with reasonable probability.

While this argument can be quite useful at a qualitative level, it
is also rather well-known that when quantitative estimates are a priority,
this method behaves rather poorly, particularly in high dimension.
One hint at this deficiency is that the assumptions do not apparently
obey any natural tensorisation property: when applied to product-form
dynamics, they do not natively produce an expected dimension-free
bound on the rate of convergence. The source of this deficiency is
not the use of drift conditions (which actually tensorise rather well),
but the use of minorisation conditions; the implied value of the constant
$\varepsilon$ degrades rather rapidly in this high-dimensional limit. 

On the basis of this observation, it was proposed to study the consequences
of replacing this minorisation condition with the `local' Poincaré
inequality on $C_{\iota}$, i.e. that for some $\kappa_{\iota}>0$,
it holds for all suitable $f$ that 
\[
\operatorname{Var}_{\pi,C_{\iota}}\left(f\right)\leq\kappa_{\iota}\pi\left(\left\Vert \nabla f\right\Vert ^{2}\right),
\]
where $\operatorname{Var}_{\pi,C_{\iota}}\left(f\right):=\inf\left\{ \pi\left(\left|f-t\right|^{2}\mathbf{1}_{C_{\iota}}\right):t\in\mathbf{R}\right\} $.
The underlying intuition then remains much the same: the drift condition
brings the dynamics into $C_{\iota}$ frequently and quickly, and
within this set, the dynamics mix rapidly.

By comparison with the usual drift-and-minorisation strategy, this
approach gives a simple proof \cite{bakry2008simple} of exponential
ergodicity with an explicit rate, typically with reasonable (if not
spectacular) estimates on the spectral gap itself. Elementary refinements
of the same proof strategy are able to verify other functional inequalities
which are each stronger than, weaker than, or incomparable to the
Poincaré inequality; see \cite{cattiaux2009lyapunov,cattiaux2010functional,cattiaux2014functional}
for details. A recent work of Taghvaei--Mehta \cite{taghvaei2021lyapunov}
helps to bring this connection full circle.

\subsection{Consequences for Mixing and Concentration}

We now move from verifying these inequalities to applying them. We
first make the convergence interpretation explicit, and then record
some less immediate consequences for covariance, concentration, and
regularisation.

\subsubsection{Marginal Convergence and Mixing Times}

Let $\rho_{t}=\mathrm{d}\mu_{t}/\mathrm{d}\pi$. If PI holds for $\pi$
with constant $\gamma$, the preceding dissipation calculation gives
\[
\chi^{2}\left(\mu_{t},\pi\right)=\operatorname{Var}_{\pi}\left(\rho_{t}\right)\leq\exp\left(-2\gamma t\right)\chi^{2}\left(\mu_{0},\pi\right).
\]
If LSI holds for $\pi$ with constant $\lambda$, then we compute
instead that
\[
\mathsf{KL}\left(\mu_{t},\pi\right)\leq\exp\left(-2\lambda t\right)\mathsf{KL}\left(\mu_{0},\pi\right).
\]
The standard comparisons 
\[
2\mathsf{TV}\left(\mu,\pi\right)^{2}\leq\mathsf{KL}\left(\mu,\pi\right),\qquad4\mathsf{TV}\left(\mu,\pi\right)^{2}\leq\chi^{2}\left(\mu,\pi\right)
\]
then allow for straightforward extraction of explicit total-variation
mixing bounds. The constant which characterises the functional inequality
determines the exponential rate of convergence, while the initial
divergence determines the cost of starting out-of-equilibrium. In
well-studied high-dimensional model problems, initial chi-squared
divergences are often controlled as $\exp\left(\mathcal{O}\left(d\right)\right)$,
whereas initial KL divergences are often of order $\mathcal{O}\left(d\right)$,
and so entropy contraction results are often capable of achieving
much tighter mixing time bounds when applicable.

For continuous-state diffusions, a point mass is singular with respect
to $\pi$, so its initial chi-squared divergence and relative entropy
are both infinite. One must therefore either assume an initial law
with finite relevant divergence or show that the dynamics make it
finite after some positive time. One useful route of this form is
provided by reverse transportation\footnote{This is to be contrasted with usual transportation-entropy inequalities,
which instead upper-bound a transportation cost by the relative entropy.}, or log-Harnack, inequalities. If $\nabla\log\pi$ is globally Lipschitz
with constant $L>0$, then by \cite{wang2010harnackboundary},
\[
\mathsf{KL}(\mu_{t},\pi)\leq\frac{L}{2(1-e^{-2Lt})}\mathcal{T}^{2}_{2}(\mu_{0},\pi),\qquad t>0,
\]
provided that the initial Euclidean quadratic transport cost $\mathcal{T}^{2}_{2}$
is finite. The coefficient is asymptotic to $1/(4t)$ as $t\downarrow0$.
In particular, a point initialisation acquires finite relative entropy
immediately when $\pi$ has finite second moment, after which an LSI
can then propagate exponential decay.

\subsubsection{Variance Bounds, Concentration Bounds, and Converse Conjectures}

Suppose that $\pi$ satisfies a Poincaré inequality with constant
$\gamma$, and write $\mathsf{C}=\operatorname{Cov}_{\pi}\left[\mathrm{Id}\right]$;
assume for simplicity that $\mathbf{E}\left[X\right]=0$. For $v\in\mathbf{R}^{d}$,
taking $G\left(x\right)=v^{\top}x$ in the Poincaré inequality yields
that $\gamma v^{\top}\mathsf{C}v\leq\left\Vert v\right\Vert ^{2}$,
and hence (by optimisation over unit $v$) that $\gamma\leq\left\Vert \mathsf{C}\right\Vert ^{-1}_{\mathrm{op}}$,
i.e. probability measures with large covariance matrices cannot have
large spectral gap. 

On the other hand, the celebrated KLS Conjecture (surveyed in \cite{alonso2015approaching})
posits the existence of a dimension-independent constant $c_{\mathrm{KLS}}\in\left(0,1\right]$
such that for all log-concave measures, one can lower-bound $\gamma\geq c_{\mathrm{KLS}}\left\Vert \mathsf{C}\right\Vert ^{-1}_{\mathrm{op}}$.
At the level of stochastic processes, this can be interpreted as saying
that in the log-concave setting, the only structural impediment to
converging to equilibrium is to traverse `long' directions in the
target. As of November 2026, a collection \cite{song20261,bizeul2026presenting,balasubramanian2026dimension}
of preprints reports claimed resolutions of the conjecture, with \cite{balasubramanian2026dimension}
in particular suggesting that one can take $c_{\mathrm{KLS}}\geq1/25$.
These works build on a string of exciting recent developments which
have substantially advanced the understanding of log-concave probability
measures in high dimension, including for example \cite{chen2021almost,klartag2023logarithmic,lee2024eldan,letwin2026kls}.

Supposing again that $\pi$ satisfies a Poincaré inequality with constant
$\gamma$, let $F:\mathbf{R}^{d}\to\mathbf{R}$ be a 1-Lipschitz function
of mean $0$, and for $\theta>0$ sufficiently small, consider $G_{\theta}:x\mapsto\exp\left(\theta F\left(x\right)\right)$.
An elementary inductive argument then shows that for $\left|\theta\right|<2\gamma^{1/2}$,
there holds the bound
\[
\mathbf{E}\left[\exp\left(\theta F\left(X\right)\right)\right]\leq\left(1-\frac{\theta^{2}}{4\gamma}\right)^{-2},
\]
whereby it follows that $\pi$ satisfies exponential concentration
for Lipschitz observables. Our earlier example exhibiting that the
exponential measure enjoys a Poincaré inequality then implies that
in general, this is the best statement of tail decay that we can expect.

Let us now seek similar conclusions under the assumption that $\pi$
satisfies a logarithmic Sobolev inequality with constant $\lambda$,
i.e. that for all suitable $G$, 
\[
\mathbf{E}\left[\frac{\left\Vert \nabla G\left(X\right)\right\Vert ^{2}}{G\left(X\right)}\right]\geq2\lambda\mathsf{Ent}\left(G\right).
\]
Considering for $v\in\mathbf{R}^{d}$ test functions of the form $G_{v}:x\mapsto\exp\left(v^{\top}x\right)$
(and assuming for simplicity that $\mathbf{E}\left[X\right]=0$),
one obtains that
\[
\lambda\log\left(\mathbf{E}\left[\exp\left(v^{\top}X\right)\right]\right)\leq\frac{\left\Vert v\right\Vert ^{2}}{2},
\]
i.e. that the moment-generating function of $X$ obeys a sub-Gaussian
bound with variance proxy proportional to $\lambda^{-1}$. By analogy
with the KLS conjecture, Bizeul \cite{bizeul2026log} formulates an
LSI analogue: 
\begin{quote}
for log-concave measures, can uniform sub-Gaussian control of linear
functionals be upgraded to a dimension-free LSI?
\end{quote}
More precisely, suppose that $X\sim\pi$ is centred and that for some
$\sigma>0$,
\[
\log\mathbf{E}\exp(v^{\top}X)\leq\frac{\sigma^{2}\|v\|^{2}}{2},\qquad v\in\mathbf{R}^{d}.
\]
The conjecture asks whether log-concavity then guarantees an LSI with
$\lambda\geq c\sigma^{-2}$ for a universal $c>0$. Bizeul \cite{bizeul2026log}
establishes the dimension-dependent bound $\lambda\geq cd^{-1/2}\sigma^{-2}$.
As of October 2026, a `raw research data' release of OpenAI \cite{openai2026subgaussian}
claims a dimension-free bound of this form.

Again under the LSI assumption, let $F:\mathbf{R}^{d}\to\mathbf{R}$
be a 1-Lipschitz function of mean $0$. A similar argument (attributed
to Ira Herbst) then establishes that for every $\theta\in\mathbf{R}$,
\[
\log\left(\mathbf{E}\left[\exp\left(\theta F\left(X\right)\right)\right]\right)\leq\frac{\theta^{2}}{2\lambda},
\]
whereby all Lipschitz functions have sub-Gaussian tails under $\pi$.
As before, we know that the Gaussian measure enjoys an LSI, and hence
no stronger tail bound can be deduced from the LSI in general.

Finally, recall the interpretation of the LSI as a contraction of
the entropy functional, i.e. 
\[
\mathsf{Ent}\left(P_{t}F\right)\leq\exp\left(-2\lambda t\right)\mathsf{Ent}\left(F\right).
\]
Somewhat remarkably, the LSI can also be shown to imply an additional
regularisation property for the semigroup, whereby the tail regularity
of the function in question is improved over time. In particular,
for $p>1$, $t>0$, define $q_{p,\lambda}\left(t\right)=1+\left(p-1\right)\exp\left(2\lambda t\right)$.
Nelson's hypercontractivity\footnote{While the term `hypercontractivity' is standard in the literature,
one should really keep in mind that it is a statement of non-expansivity,
i.e. `contractivity' with pre-factor $1$. Indeed, by considering
$F\equiv1$, one sees that strict contractivity cannot hold for all
functions.} theorem then states that 
\[
\left\Vert P_{t}F\right\Vert _{q_{p,\lambda}\left(t\right)}\leq\left\Vert F\right\Vert _{p}.
\]
That is, if the function $F$ is initially $\mathrm{L}^{p}$-integrable
under $\pi$, then as $t$ grows, the function $P_{t}F$ becomes $\mathrm{L}^{q}$-integrable
for each $q>p$ in finite time (and even rather quickly!). This reflects
that under LSI the semigroup regularises as well as relaxes: it improves
integrability while moving functions towards equilibrium. Moreover,
it can be shown that at the right level of generality, this hypercontractivity
property implies the LSI, i.e. the two properties are equivalent.
This is detailed neatly in Chapter 2 of \cite{chafai2023logarithmic},
and is even explained rather clearly in the early work of Gross \cite{gross1975logarithmic}
which gave the LSI its name.

\subsection{Horizons: Beyond Poincaré and Log-Sobolev}

The preceding discussion has focused predominantly on the Poincaré
and log-Sobolev inequalities, due to their relative utility and applicability.
The following is intended to provide some informative (if brief) glimpses
at the broader literature on energy-entropy inequalities, explaining
how they enrich the basic picture in terms of scope and consequences.
A first reading of this article may safely proceed directly to the
conclusion.

\paragraph{$\mathrm{L}^{1}$-type Inequalities}

Recall the Poincaré inequality, which states that for mean-centred
functions $F$, the $\mathrm{L}^{2}$ norm of $F$ can be controlled
by the $\mathrm{L}^{2}$ norm of $\nabla F$. Cheeger-type inequalities
instead impose that for \emph{median-centred} $F$, one can control
\[
\mathbf{E}\left[\left|F\right|\right]\leq C\mathbf{E}\left[\left\Vert \nabla F\right\Vert \right].
\]
By carefully considering suitable functions of the form $F=G^{2}$,
these can be seen to imply Poincaré inequalities (with different but
explicit constants); moreover, under log-concavity, the equivalence
flows both ways.

Structurally, a key aspect of Cheeger inequalities is that by working
with more permissive $\mathrm{L}^{1}$-type norms, one can work with
`rougher' test functions $F$ and obtain finer insights about the
geometry of $\pi$. In particular, for reasonable sets $A$, one can
construct sequences of differentiable approximations $F$ to the function
$\mathbf{1}_{A}$ to which the Cheeger inequality can be applied,
and obtain comparisons of the form
\[
\min\left\{ \pi\left(A\right),\pi\left(A^{\complement}\right)\right\} \leq C\pi^{+}\left(A\right),
\]
where $\pi^{+}$ denotes a suitable notion of boundary measure (the
`Minkowski content' of $A$ with respect to $\pi$). In this respect,
Cheeger inequalities are closely linked to the isoperimetric problem
for probability measures.

Further in this direction, `weak' Cheeger inequalities \cite{bobkov2007large,cattiaux2010functional}
take the form
\[
\mathbf{E}\left[\left|F\right|\right]\leq\alpha\mathbf{E}\left[\left\Vert \nabla F\right\Vert \right]+\beta\mathrm{Osc}\left(F\right),\qquad\left(\alpha,\beta\right)\in\mathcal{S}_{\mathrm{WCI}}
\]
where $\mathrm{Osc}\left(F\right)=\sup F-\inf F$, $\sup$ and $\inf$
are computed essentially under $\pi$, and $\mathcal{S}_{\mathrm{WCI}}\subseteq\mathbf{R}^{2}_{+}$
denotes some `Pareto frontier' along which one can trade off the two
quantities on the right-hand side. By contrast with both the Poincaré
and Cheeger inequalities, these inequalities can hold even for heavy-tailed
measures, and can be used to take a functional perspective on the
analysis of their geometric and isoperimetric properties \cite{bobkov1997some,bobkov1999isoperimetric}.

\paragraph{Weakened and Strengthened $\mathrm{L}^{2}$-type Inequalities}

Recall the Poincaré inequality, which, as discussed previously, entails
exponential concentration for Lipschitz observables. Naturally, this
precludes its application to more heavy-tailed problems.

The `weak' Poincaré inequality \cite{rockner2001weak,wang2002coupling,cattiaux2007weak,andrieu2026weak}
follows the weak Cheeger inequality in adding a defect to this inequality,
taking the form
\[
\mathbf{E}\left[\left|F\right|^{2}\right]\leq\alpha\mathbf{E}\left[\left\Vert \nabla F\right\Vert ^{2}\right]+\beta\mathrm{Osc}\left(F\right)^{2},\qquad\left(\alpha,\beta\right)\in\mathcal{S}_{\mathrm{WPI}}
\]
again for all mean-centred $F$ and some target-dependent family of
admissible pairs $\mathcal{S}_{\mathrm{WPI}}$. Such inequalities
are available for many heavy-tailed measures for which PI fails, and
can be used to prove slower-than-exponential rates of convergence
for the corresponding Langevin diffusion in the $\chi^{2}$ divergence,
where the precise rate depends on the trade-off described by $\mathcal{S}_{\mathrm{WPI}}$.

On the other hand, the `super' Poincaré inequality \cite{wang2000functional,cattiaux2009lyapunov,zitt2011super}
adds an additional term to the right-hand side which is much \emph{smaller}
than the second moment, instead taking the form
\[
\mathbf{E}\left[\left|F\right|^{2}\right]\leq\alpha\mathbf{E}\left[\left\Vert \nabla F\right\Vert ^{2}\right]+\beta\mathbf{E}\left[\left|F\right|\right]^{2},\qquad\left(\alpha,\beta\right)\in\mathcal{S}_{\mathrm{SPI}},
\]
now for general $F$ of suitable regularity. Taking $F\equiv1$, one
sees that $\beta\geq1$. Supposing that there exists $\left(\alpha,\beta\right)\in\mathcal{S}_{\mathrm{SPI}}$
with $\beta<2$, one can infer a classical Poincaré inequality for
$\pi$. Indeed, for mean-centred $G$, set $G_{+}=\max\{G,0\}$ and
$G_{-}=\max\{-G,0\}$. Applying SPI separately to these functions,
using the chain rule, and checking that $\mathbf{E}\left[G_{+}\right]=\mathbf{E}\left[G_{-}\right]=\tfrac{1}{2}\mathbf{E}\left[\left|G\right|\right]$,
one can argue that
\[
\mathbf{E}\left[G^{2}\right]\leq\alpha\mathbf{E}\left[\left\Vert \nabla G\right\Vert ^{2}\right]+\frac{\beta}{2}\left(\mathbf{E}\left[\left|G\right|\right]\right)^{2}\leq\alpha\mathbf{E}\left[\left\Vert \nabla G\right\Vert ^{2}\right]+\frac{\beta}{2}\mathbf{E}\left[G^{2}\right].
\]
Rearranging then yields a Poincaré inequality with $\gamma\geq(2-\beta)/(2\alpha)$.

Pursuing this perspective further, one can show that for suitable
frontiers $\mathcal{S}_{\mathrm{SPI}}$, these estimates can even
imply sharper estimates of mixing time than PI alone provides, often
demonstrating an anomalous `faster-than-exponential' rate of convergence
when far from equilibrium, before later slowing down to the usual
exponential rate of relaxation. In effect, by allowing for larger
values of $\beta$, one can take $\alpha$ much smaller, yielding
`effective contraction at rate $\alpha^{-1}$', up until a certain
point. We note also a connection to other geometric notions including
Faber-Krahn inequalities and the `Spectral Profile' of Grigor'yan
and others; see e.g. \cite{grigor1994heat,grigor2006heat,goel2006mixing,grigoryan2009heat}
for more on this perspective.

\paragraph{Weakened and Strengthened Entropic Inequalities}

Recall the logarithmic Sobolev inequality. The `weak' LSI \cite{cattiaux2007weak}
expands this to the defective form
\[
\mathsf{Ent}\left(F^{2}\right)\leq\alpha\mathbf{E}\left[\left\Vert \nabla F\right\Vert ^{2}\right]+\beta\mathrm{Osc}\left(F\right)^{2},\qquad\left(\alpha,\beta\right)\in\mathcal{S}_{\mathrm{WLSI}},
\]
which allows for the quantitative study of convergence to equilibrium
in relative entropy, even for target distributions with heavier-than-Gaussian
tails. 

Conversely, the `super' LSI \cite{gross1998herbst} refines this to
\[
\mathsf{Ent}\left(F^{2}\right)\leq\alpha\mathbf{E}\left[\left\Vert \nabla F\right\Vert ^{2}\right]+\beta\mathbf{E}\left[\left|F\right|^{2}\right],\qquad\left(\alpha,\beta\right)\in\mathcal{S}_{\mathrm{SLSI}},
\]
and is closely connected with strong short-time regularisation properties
of the semigroup, i.e. reinforcements of the hypercontractivity property
enjoyed by the LSI.

\paragraph{Weighted Sobolev-Type Inequalities}

Recall the Poincaré inequality. In `weighted' Poincaré inequalities
\cite{bobkov2009weighted,bobkov2009weighted2,saumard2019weighted,huguet2024poincare},
one seeks a positive weight $\omega:\mathbf{R}^{d}\to\mathbf{R}_{+}$
such that for mean-centred $F$, it rather holds that
\[
\mathbf{E}\left[\left|F\left(X\right)\right|^{2}\right]\leq\mathbf{E}\left[\omega\left(X\right)\left\Vert \nabla F\left(X\right)\right\Vert ^{2}\right].
\]
In practice, such inequalities are typically most interesting when
the basic Poincaré inequality fails, and so one has to take $\omega$
unbounded. It is then of particular interest to identify the `minimal'
such weight, i.e. the gentlest perturbation which allows for the inequality
to hold.

Dynamically, this energy corresponds to a $\pi$-reversible diffusion
with diffusion coefficient $\sqrt{2\omega(x)}$. When $\mathbf{E}_{\pi}\omega<\infty$,
it can also be viewed as a time-change of the ordinary Langevin diffusion
with invariant probability density $\pi^{\omega}\propto\pi\omega$.
The resulting generator and Dirichlet form for this process are then
given as
\begin{align*}
\mathcal{L}_{\omega}F & =\pi^{-1}\nabla\cdot\left(\pi\omega\nabla F\right)\\
 & =\left\langle \omega\nabla\log\pi+\nabla\omega,\nabla F\right\rangle +\omega\Delta F,\\
\mathcal{E}_{\omega}\left(F,G\right) & =\mathbf{E}\left[\omega\left(X\right)\left\langle \nabla F\left(X\right),\nabla G\left(X\right)\right\rangle \right].
\end{align*}
Under the usual regularity and non-explosion assumptions, the weighted
Poincaré inequality $\mathcal{E}_{\omega}\left(F,F\right)\geq\gamma_{\omega}\mathrm{Var}_{\pi}\left(F\right)$
entails that this process relaxes to the equilibrium distribution
$\pi$ in chi-squared divergence with exponential rate $\geq2\gamma_{\omega}$.

Such perspectives have been used recently in the development of improved
algorithms for stochastic simulation; see e.g. \cite{bertazzi2025sampling,mijatovic2026tempered}.
As mentioned above, focus has largely centred around the case of unbounded
weights with applications to heavy-tailed sampling, for which well-designed
weights can induce exponential convergence in situations where the
basic Langevin diffusion would not have the same favourable properties.

Take as an example the (rescaled) univariate Student distribution
with $\nu>0$ degrees of freedom, with density $\pi_{\nu}\left(x\right)\propto\left(1+x^{2}\right)^{-\left(\nu+1\right)/2}$
on $\mathbf{R}$; for simplicity, focus on the case $\nu\geq3$. The
polynomial tails mean that no nontrivial Poincaré inequality can hold
for $\pi_{\nu}$. Nevertheless, Proposition 7 of \cite{mousavi2023towards},
building on Bobkov-Ledoux \cite{bobkov2009weighted2}, shows that
a weak Poincaré inequality holds for $\pi_{\nu}$, where $\mathcal{S}_{\mathrm{WPI}}$
includes admissible pairs
\[
\left(\alpha,\beta\right)=\left(\frac{2}{\nu}+2\left(1+\frac{1}{\nu}\right)s^{-2},s^{\nu}\right),\qquad s>0.
\]
For initialisations $\mu$ satisfying $\left\Vert \mathrm{d}\mu/\mathrm{d}\pi_{\nu}\right\Vert _{\infty}<\infty$,
routine calculations (following e.g. \cite{rockner2001weak}) then
yield a trend to equilibrium of $\chi^{2}\left(\mu_{t,\nu},\pi_{\nu}\right)\in\mathcal{O}\left(t^{-\nu/2}\right)$.

On the other hand, the weighted Langevin diffusion with weight $\omega\left(x\right)=1+x^{2}$
fares rather better. The dynamics take the form
\[
\mathrm{d}X_{t}=-\left(\nu-1\right)X_{t}\,\mathrm{d}t+\sqrt{2\left(1+X^{2}_{t}\right)}\,\mathrm{d}W_{t},
\]
with a strong restoring drift pushing trajectories back towards the
origin. Theorem 3.1 of Bonnefont-Joulin-Ma \cite{bonnefont2016note}
and Section 5 of Huguet \cite{huguet2024poincare} then give the (sharp)
weighted inequality below for $\pi_{\nu}$-centred $F$; see also
Corollary 14 of Nguyen \cite{nguyen2013dimensional} for the parameter
range considered here. 
\[
\left(\nu-1\right)\mathbf{E}_{\nu}\left[\left|F\left(X\right)\right|^{2}\right]\leq\mathbf{E}_{\nu}\left[\left(1+X^{2}\right)\left\Vert \nabla F\left(X\right)\right\Vert ^{2}\right].
\]
Writing $\mu^{\omega}_{t,\nu}$ for the law of this weighted diffusion,
we therefore obtain
\[
\chi^{2}(\mu^{\omega}_{t,\nu},\pi_{\nu})\leq e^{-2(\nu-1)t}\chi^{2}(\mu_{0},\pi_{\nu}),
\]
provided the initial chi-squared divergence is finite.

Weighted Poincaré inequalities also provide an input to hypocoercive
convergence estimates for non-reversible lifts; see Mijatović-Niraj
\cite{mijatovic2026hypocoercivity} for tempered bouncy particle samplers
with heavy-tailed targets.

One can equally imagine similar statements concerning dynamics involving
non-scalar noise, which would then be understood as `Riemannian' Poincaré
inequalities with respect to a suitable metric; these are yet more
complicated, but also insightful in certain contexts \cite{kolesnikov2016riemannian,chewi2020exponential,cui2024optimal}.
Equally, one can seek weighted variants of the LSI, admitting similar
interpretations and applications \cite{cattiaux2011some,huguet2026logarithmic}.

\paragraph{Łojasiewicz-type Inequalities in Continuous Optimisation}

Returning to our `reverse-Grönwall' interpretation of functional inequalities,
one can even return all the way to conventional optimisation on smooth
spaces by gradient flow methods. Given a smooth and lower-bounded
$F:\mathbf{R}^{d}\to\mathbf{R}$ (without loss of generality, one
may assume by shifting that $\inf F=0$), the gradient flow dynamics
for minimising $F$ take the form
\[
\dot{X}_{t}=-\nabla F\left(X_{t}\right),
\]
and so if it holds globally that $F$ is `gradient-dominated' in the
sense that for some $m>0$ and for all $x\in\mathbf{R}^{d}$
\[
\left\Vert \nabla F\left(x\right)\right\Vert ^{2}\geq2mF\left(x\right),
\]
then one deduces by the usual argument that
\[
F\left(X_{t}\right)\leq\exp\left(-2mt\right)F\left(X_{0}\right),
\]
i.e. exponential convergence of the objective gap along the flow.
Such conditions are now commonly called Polyak--Łojasiewicz inequalities
(see \cite{karimi2016linear}), and are widely used as a tractable
relaxation of strong convexity for which gradient methods nevertheless
converge reliably. 

The precise relation of these conditions to classical convexity is
subtle and remains an active topic of research. Although weaker than
strong convexity, recent work has revealed that when combined with
additional regularity assumptions, global PŁ inequalities impose substantial
geometric restrictions on the function. For instance, \cite{nejma2025polyak}
shows that a $C^{2}$ function satisfying a global PŁ inequality and
having a set of minimisers which is bounded and nonempty necessarily
has a \emph{unique} minimiser and is \emph{strongly convex} on a sufficiently
low sublevel set. See also \cite{boumal2026smooth,sontag2026smooth}
for some related recent works which seek to shed similar light on
the structure of the class of PŁ functions.

Finally, by analogy with the `weak' inequalities of previous sections,
note that any comparison of the form
\[
\left\Vert \nabla F\left(x\right)\right\Vert ^{2}\geq\Psi\left(F\left(x\right)\right)
\]
with $\Psi\left(r\right)>0$ for $r>0$ can furnish some (perhaps
slower-than-exponential) convergence rates for the gradient flow dynamics;
such inequalities are typically known as Kurdyka--Łojasiewicz inequalities
and hold locally for broad classes of definable functions; see e.g.
\cite{bolte2010characterizations} for the general picture. 

A concrete example arises in Monmarché’s study of the continuous Curie--Weiss
model \cite{monmarche2025free}, a toy model of statistical mechanics
which exhibits phase transitions as a certain temperature parameter
varies. The mean-field dynamics associated to this system can be interpreted
as a Wasserstein gradient flow of a free energy functional $\mathcal{F}$
on the space of probability measures. Above the critical temperature,
$\mathcal{F}$ enjoys a PŁ inequality, whereby the associated dynamics
are known to converge to their equilibrium at an exponential rate,
per Remark 9 of \cite{monmarche2024local}. On the other hand, at
criticality, the PŁ inequality (with any positive constant) fails,
but a KŁ inequality with $\Psi\left(r\right)\asymp\min\left\{ r^{3/2},r\right\} $
persists, yielding decay of the excess free energy at a polynomial-in-time
rate. In particular, the use of generalised coercivity inequalities
of this type is able to provide non-trivial guarantees on the convergence
behaviour of optimisation dynamics in rather degenerate scenarios.

\paragraph{Functional Inequalities beyond `Energy-Entropy' Comparisons}

The discussion here has focused primarily on inequalities of `energy-entropy'
type, i.e. comparing the nominal objective function (`entropy') to
its infinitesimal dissipation (`energy'), as these are of primary
interest for studying the long-time properties of Markov processes.
Nevertheless, it bears mentioning that from the perspective of geometric
functional inequalities at large, these are just one rung in a more
expansive ladder of inequalities, cutting across curvature, isoperimetry,
transportation-entropy inequalities, concentration inequalities, and
more; see Section 3.5 of \cite{bogachev2012monge} for a detailed
review of this `Bogachev-Kolesnikov hierarchy'. 

In general, there is a clearly directed hierarchy of implications
between these different families of inequality, passing strictly in
one direction. Under the assumption of log-concavity for the probability
measure in question, stunning results of E. Milman \cite{milman2009concentration,milman2009role,milman2012properties}
show that several of these one-way implications can be reversed up
to universal constants. Concentration, isoperimetry, and functional
inequalities of corresponding strength thereby become quantitatively
comparable when specialised to this setting, despite being genuinely
distinct for arbitrary measures.

\section{Literature Guide}

The literature around these topics is vast, with contributions from
researchers coming from probability, functional analysis, optimal
transport, and much besides. The most useful entry point depends on
what the reader already knows and what they hope to do next. The suggestions
below are organised accordingly; none is intended to be a prerequisite
for all the others.

\paragraph{For a first course in semigroups and curvature.}

The textbook of Bakry--Gentil--Ledoux \cite{bakry2014analysis}
is the standard comprehensive reference for diffusion semigroups,
the $\Gamma$-calculus, curvature, and the consequences of PI and
LSI. The notes of Michel Bonnefont \cite{bonnefont2022poincare} are
shorter and more targeted, developing a varied toolkit for quantitative
Poincaré inequalities. The notes of Chafaï--Lehec \cite{chafai2023logarithmic}
give a focused route into logarithmic Sobolev inequalities and their
variants, including their incarnations in discrete spaces.

\paragraph{For high-dimensional probability, transport, and sampling.}

Ramon van Handel's lecture notes \cite{van2016probability} place
functional inequalities inside a broader course on probability in
high dimension and are particularly effective at conveying the fundamental
interest of seeking dimension-free results. Xiaohui Chen's notes \cite{chen2020geometric}
approach functional inequalities through PDE methods, and emphasise
connections with optimal transport. Sinho Chewi's authoritative text
on log-concave sampling \cite{chewi2026logconcave} develops many
of these ideas with a firmer focus on their applications to non-asymptotic
analysis of sampling algorithms.

\paragraph{For spectral theory, Markov chains on discrete spaces, and mixing
times.}

The books of Mu-Fa Chen \cite{chen2004markov,chen2005eigenvalues}
situate functional inequalities in a classical spectral and functional-analytic
context. Montenegro--Tetali \cite{montenegro2006mathematical} write
for an audience interested in finite-state space chains and theoretical
computer science, with particular attention to mixing time estimates
from arbitrary initialisation. The book of Feng-Yu Wang \cite{wang2005functional}
is broader, connecting many functional inequalities with Markov semigroups
and stochastic processes on manifolds. Among other things, these references
help to clarify which aspects of the story told in this article are
able to adapt and survive beyond the classical setting of diffusions
on Euclidean spaces.

\paragraph{For coupling and Lyapunov methods.}

Andreas Eberle's lecture notes \cite{eberle2026markov} provide a
broad probabilistic entry point linking semigroups, mixing, sampling,
couplings, and functional inequalities. Yann Ollivier's influential
work \cite{ollivier2009ricci} on coarse Ricci curvature is an accessible
entry into the idea that metric contraction can imply functional inequalities,
as well as a rich range of other consequences. The surrounding literature
provides a fertile ground for understanding the impact of these results
by example; see e.g. \cite{chen1994application,chen1998trilogy,eberle2016reflection,salez2025modern}
for some examples. Readers arriving from Foster--Lyapunov theory
would be well served by beginning with the initial work of Bakry--Barthe--Cattiaux--Guillin
on deriving Poincaré inequalities from drift and local estimates \cite{bakry2008simple},
followed by the extensions and survey chapter of Cattiaux--Guillin
\cite{bakry2008rate,cattiaux2009lyapunov,cattiaux2010functional,cattiaux2014functional}.
This route is closest in spirit to classical stochastic stability
and makes the complementarity advocated in the prelude especially
concrete.

\paragraph{For heavy tails and the wider family of inequalities.}

The works of Röckner--Wang \cite{rockner2001weak} and Cattiaux--Gozlan--Guillin--Roberto
\cite{cattiaux2010functional} provide routes into weak, super, and
weighted inequalities. These provide context on why PI and LSI serve
as natural boundaries for many applications, and on how these modified
inequalities become necessary and useful beyond these boundaries.
The literature connecting functional inequalities to isoperimetry,
transport, and concentration is still larger; the works of Emanuel
Milman \cite{milman2009concentration,milman2009role,milman2012properties},
Sergei Bobkov \cite{bobkov1999isoperimetric,bobkov2005entropy}, and
Michel Ledoux \cite{ledoux2011concentration} offer excellent entry
points.

\section{Denouement}

Functional inequalities describe long-time convergence at a deliberately
coarse scale, replacing a detailed account of trajectories by a comparison
between a functional and the rate at which it dissipates over time.
For the Langevin diffusion, this turns the Poincaré inequality and
the logarithmic Sobolev inequality into transparent statements about
contraction of chi-squared divergence and relative entropy respectively.
Beyond these questions of convergence, the same viewpoint is rich
enough to also shed light on their consequences for concentration,
regularisation, and spectral properties.

Operationally, a good functional inequality is a robust way to say
that
\begin{quote}
\begin{center}
this Markov process mixes well at large scales,
\par\end{center}

\end{quote}
without committing to a complete local description of the dynamics.
This is why the inequalities compose so effectively: the combination
of measures by taking products, transport, convolution, bounded perturbations,
and much more can be analysed seamlessly without rebuilding a trajectory-level
argument from the beginning.

The cost of this coarse-graining is that the first useful inequality
can be hard to obtain. Uniform curvature gives a remarkably clean
answer but is often too rigid as an assumption, transfer principles
require a suitable reference model, and coupling or Lyapunov arguments
still demand genuine probabilistic insight; there is no unqualified
escape from the principle of `conservation of mathematical difficulty'.
The gain is that once this initial work has been done, the conclusion
can be propagated and reused with unusual efficiency.

My hope is that this perspective makes the subject feel less like
a collection of analytic incantations associated to the names of eminent
mathematical analysts of yesteryear, and more like a cohesive and
worthwhile language for studying the large-scale and long-time properties
of Markov processes. Functional inequalities are a flexible complement
to pathwise analysis: they interface meaningfully with classical local
and probabilistic arguments, enable the conversion of these conclusions
into global information, and enjoy stability to a range of perturbations
which might be more challenging to study directly by probabilistic
means. In these regards, they offer great value as a tool in the modern
toolbox of Markov chain analysis.

\section*{Acknowledgements}

Much of what I learned about these topics was catalysed and supported
by a rewarding extended collaboration with Christophe Andrieu, Anthony
Lee, and Andi Wang. Conversations and correspondence with numerous
other experts on functional inequalities have also been instructive
in shaping my broader understanding of the area.

I am grateful to the organisers of the INI programme ``Stochastic
Systems for Anomalous Diffusion'' (in particular, Codina Cotar, Aleks
Mijatovic, and Andrew Wade) for their encouragement to deliver a mini-course
on this topic, and to write this article as a supplement. I also gratefully
acknowledge the CSML Reading Group in the University of Warwick Department
of Statistics, and the participants of the 2025 MATRIX Institute workshop
on ``Bayesian Learning of Very-High-Dimensional Physical Process
Models'' for providing opportunities to deliver additional iterations
of the mini-course and refine its presentation.

{\footnotesize\bibliographystyle{abbrv}
\bibliography{power_functional_inequalities}
}{\footnotesize\par}
\end{document}